\documentclass[fleqn,twoside,10pt,a4paper]{article}

\usepackage{amssymb}
\usepackage{epsfig}

\newtheorem{theorem}{Theorem}
\newtheorem{lemma}[theorem]{Lemma}

\title{
  An Infinite Suite of Links--Gould Invariants
}

\author{
  David~~De Wit%
    \footnote{
    RIMS, Kyoto University 606-8502, Japan.
    \texttt{ddw@kurims.kyoto-u.ac.jp}
  }
}

\begin{document}

\maketitle

\begin{abstract}
  \noindent
  This paper describes a method to obtain state model parameters for an
  infinite series of Links--Gould link invariants $LG^{m,n}$, based on
  quantum R matrices associated with the
  $(\dot{0}_m\,|\,\dot{\alpha}_n)$ representations of the quantum
  superalgebras $U_q[gl(m|n)]$. Explicit details of the state models
  for the cases $n=1$ and $m=1,2,3,4$ are supplied.  Some gross
  properties of the link invariants are provided, as well as some
  explicit evaluations.
\end{abstract}


\section{Overview}

In 1992, Jon Links and Mark Gould \cite{LinksGould:92b} described a
method for constructing link invariants from quantum superalgebras. That
work stopped short of evaluations of the invariants due to want of an
efficient computational method.  In 1999, the author, in collaboration
with Jon Links and Louis Kauffman \cite{DeWitKauffmanLinks:99a}, first
evaluated a two-variable example of one these invariants, using a state
model.  We used the $(0,0\,|\,\alpha)$ representations of
$U_q[gl(2|1)]$, and labeled our resulting $(1,1)$-tangle invariant
$LG$, `\emph{the} Links--Gould invariant'. In that paper, and
subsequently in \cite{DeWit:99a}, we showed that whilst $LG$ would
detect neither inversion nor mutation, it was still able to distinguish
all prime knots of up to $10$ crossings, making it more powerful than
the HOMFLY and Kauffman invariants.

Here, we generalise the notation, denoting $LG^{m,n}$ as ``the
Links--Gould invariant associated with the
$(\dot{0}_m\,|\,\dot{\alpha}_n)$ representation of $U_q[gl(m|n)]$''.
For the case $n=1$, we will write $LG^{m}\equiv LG^{m,1}$, so our
previous invariant $LG$ was in fact $LG^{2}$.  This generalisation is
motivated by the automation of a procedure to construct the appropriate
R matrices \cite{DeWit:99c,DeWit:99d}; previously, we were limited to
the $m=2$ case, for which the R matrix had been calculated by hand.

We explicitly demonstrate the construction of state model parameters
for $LG^{m,n}$, illustrating our results for $LG^m$, for the cases
$m=1,2,3,4$. Further, we describe some of the gross properties of these
invariants, and provide a limited set of evaluations of them.

Although these invariants $LG^{m,n}$ are not more powerful in their
gross properties than $LG^2$ (they can detect neither inversion nor
mutation), each one is expected to distinguish many more knots $K$ as
the degree of the polynomials $LG^{m,n}_K$ increases rapidly with $m$
and $n$.  Perhaps more significantly, the development of the current
formalism points the way towards automation of the evaluation of more
general classes of quantum link invariants; a discussion of this is
provided.


\section{Quantum superalgebra state models}

Corresponding to each finite dimensional highest weight representation
of each quantum superalgebra, there exists a quantum link invariant
($LG$), originally described in \cite{LinksGould:92b}. These invariants
are similar to those associated with the usual (i.e. ungraded) quantum
algebras (e.g. \cite{Kauffman:87b,Kauffman:93,ReshetikhinTuraev:90}),
although there are some technical differences.

Here, we describe the construction of parameters for state models for
evaluating a class of these invariants. Specifically, we will define
$LG^{m,n}$ to be the quantum link invariant associated with the
representation $\pi\equiv\pi_\Lambda$ of highest weight
$\Lambda=(\dot{0}_m\,|\,\dot{\alpha}_n)$, of the quantum superalgebra
$U_q[gl(m|n)]$.  To do so, we first broadly introduce the algebraic
structures, then we briefly review the terminology used to describe
state model parameters, and finally, we look at the construction of
specific state model parameters for our particular class of
representations.


\subsection{The quantum superalgebra $U_q[gl(m|n)]$}
\label{sec:quantumsuperalgebra}

$U_q[gl(m|n)]$ is a unital super (i.e. $\mathbb{Z}_2$-graded) algebra
with free parameter $q$. In the limit $q\to 1$, it degenerates to the
ordinary Lie superalgebra $gl(m|n)$.  Here, we provide a broad outline
of $U_q[gl(m|n)]$ in terms of generators and relations, for readers not
familiar with it. This material is largely abstracted from the fuller
description contained in \cite{Zhang:93} (see also \cite{DeWit:99c}).


\subsubsection{$U_q[gl(m|n)]$ generators}

A set of generators for $U_q[gl(m|n)]$ is:
\begin{eqnarray*}
  \left\{
  \! \!
    \begin{array}{rll}
      K_{a},       & 1 \leqslant a \leqslant m+n & \mathrm{Cartan}
        \\[0.5mm]
      {E^{a}}_{b}, & 1 \leqslant a < b \leqslant m+n & \mathrm{raising}
        \\[0.5mm]
      {E^{b}}_{a}, & 1 \leqslant a < b \leqslant m+n & \mathrm{lowering}
    \end{array}
  \! \!
  \right\}.
\end{eqnarray*}
An equivalent notation for $K_{a}$ is $q_a^{{E^a}_a}$, where we have
introduced the notation $q_a\triangleq q^{{(-)}^{[a]}}$.  For any power
$N$, we may write $q_a^N$, and hence $K_{a}^N$, thus for $ M,N \in
\mathbb{C}$:
\begin{eqnarray*}
  K_{a}^{M}
  K_{a}^{N}
  =
  K_{a}^{M+N}
  \quad
  \mathrm{where}
  \quad
  K_{a}^0
  \equiv
  \mathrm{Id},
  \qquad
\end{eqnarray*}
where $\mathrm{Id}$ is the $U_q[gl(m|n)]$ identity element.

Using the following $\mathbb{Z}_2$ grading on the $gl(m|n)$
\emph{indices}:
\begin{eqnarray*}
  [a]
  \triangleq
  \left\{
  \begin{array}{lll}
    0 \quad & \mathrm{if} \quad 1 \leqslant a \leqslant m \qquad\qquad
            & \mathrm{even}
    \\
    1       & \mathrm{if} \quad m+1 \leqslant a \leqslant m+n    &
            \mathrm{odd},
  \end{array}
  \right.
\end{eqnarray*}
we may define a natural $\mathbb{Z}_2$ grading on the generators:
\begin{eqnarray*}
  [K_{a}^N]
  \triangleq
  0,
  \qquad \qquad
  [{E^a}_b]
  \triangleq
  [a] + [b]
  \quad
  (\mathrm{mod}\;2),
\end{eqnarray*}
and we use the terms ``even'' and ``odd'' for generators in the same
manner as we do for indices.  Elements of $U_q[gl(m|n)]$ are said to be
\emph{homogeneous} if they are linear combinations of generators of the
same grading.  The product $XY$ of homogeneous $X,Y\in U_q[gl(m|n)]$
has grading:
\begin{eqnarray*}
  [ X Y ]
  \triangleq
  [X] + [Y]
  \quad
  ( \mathrm{mod}\;2 ).
\end{eqnarray*}
Within the full set of generators, we have the
$U_q[gl(m|n)]$ \emph{simple} generators:
\begin{eqnarray*}
  \left\{
    \! \! \!
    \begin{array}{rll}
      K_{a},         & 1 \leqslant a \leqslant m+n & \mathrm{Cartan} \\
      {E^{a}}_{a+1}, & 1 \leqslant a < m+n & \mathrm{simple~raising} \\
      {E^{a+1}}_{a}, & 1 \leqslant a < m+n & \mathrm{simple~lowering}
    \end{array}
    \! \! \!
  \right\},
\end{eqnarray*}
such that the remaining nonsimple generators may be expressed in terms
of these \cite[p1238,~(2)]{Zhang:93}.  The fact that there are $m+n-1$
simple raising generators indicates that $U_q[gl(m|n)]$ has rank
$m+n-1$.


\subsubsection{$U_q[gl(m|n)]$ relations}

The \emph{graded commutator}
$
  [\cdot,\cdot]
  :
  U_q[gl(m|n)] \times U_q[gl(m|n)]
  \to
  U_q[gl(m|n)]
$,
is defined for homogeneous
$ X, Y \in U_q[gl(m|n)] $ by:
\begin{eqnarray*}
  [ X, Y ]
  \triangleq
  X Y - {(-)}^{[X][Y]} Y X,
\end{eqnarray*}
and extended by linearity.  With this, we have the following
$U_q[gl(m|n)]$ relations:

\begin{enumerate}
\item
  The Cartan generators all commute:
  \begin{eqnarray*}
    K_{a}^{M}
    K_{b}^{N}
    =
    K_{b}^{N}
    K_{a}^{M},
    \qquad \qquad
    M, N \in \mathbb{C}.
  \end{eqnarray*}

\item
  The Cartan generators commute with the simple raising and lowering
  generators in the following manner:
  \begin{eqnarray*}
    K_{a}
    {E^b}_{b\pm1}
    =
    q_a^{(\delta^a_b - \delta^a_{b\pm1} )}
    {E^b}_{b\pm1}
    K_{a}.
  \end{eqnarray*}

\item
  The squares of the odd simple generators are zero:
  \begin{eqnarray*}
    {( {E^{m}}_{m+1} )}^2
    =
    {( {E^{m+1}}_{m} )}^2
    =
    0.
  \end{eqnarray*}
  (This implies that the squares of nonsimple odd generators are
  also zero.)

\item
  The non-Cartan generators satisfy the following commutation relations:
  \begin{eqnarray*}
    [ {E^a}_{a+1}, {E^{b+1}}_b ]
    =
    \delta^a_b
    \frac{
      K_{a} K_{a+1}^{-1} - K_{a}^{-1} K_{a+1}
    }{
      q_a - \overline{q}_a
    },
  \end{eqnarray*}
  where we have written $\overline{q}\equiv q^{-1}$ for brevity.
  We also have, for $|a-b| > 1$, the commutations:
  \begin{eqnarray*}
    \hspace{-27pt}
    {E^{a}}_{a+1}
    {E^{b}}_{b+1}
    =
    {E^{b}}_{b+1}
    {E^{a}}_{a+1}
    \qquad
    \mathrm{and}
    \qquad
    {E^{a+1}}_{a}
    {E^{b+1}}_{b}
    =
    {E^{b+1}}_{b}
    {E^{a+1}}_{a}.
  \end{eqnarray*}

\item
  Lastly, we have the $U_q[gl(m|n)]$ \emph{Serre relations}; their
  inclusion ensures that the algebra is reduced enough to be
  \emph{simple}. We omit these for brevity; they are not required below.

\end{enumerate}


\subsubsection{$U_q[gl(m|n)]$ as a Hopf superalgebra}

When equipped with an appropriate%
\footnote{%
  The details of these structures are not required here; the
  reader can find them in \cite{DeWit:98,DeWit:99c}.
}
coproduct $\Delta$, counit $\varepsilon$ and antipode $S$, we may
regard $U_q[gl(m|n)]$ as a quasitriangular Hopf superalgebra.
This means that it possesses an R matrix $\check{R}$, an operator on
the tensor product $U_q[gl(m|n)]\otimes U_q[gl(m|n)]$, satisfying
the quantum Yang--Baxter equation (QYBE) in the form:
\begin{equation}
  (\check{R} \otimes I)
  (I \otimes \check{R})
  (\check{R} \otimes I)
  =
  (I \otimes \check{R})
  (\check{R} \otimes I)
  (I \otimes \check{R}),
  \label{eq:QYBEinBraidForm}
\end{equation}
immediately recognisable as the braid relation:
\begin{equation}
  \sigma_1
  \sigma_2
  \sigma_1
  =
  \sigma_2
  \sigma_1
  \sigma_2.
  \label{eq:BraidRelation}
\end{equation}


\subsection{State model parameters}
\label{sec:Statemodelparameters}

The following comments briefly describe in what is well-documented in
the literature (see, e.g.  \cite{Hennings:91,Kauffman:97a}). They are
included so as to introduce our particular notation.

\begin{figure}[htbp]
  \begin{center}
    \input{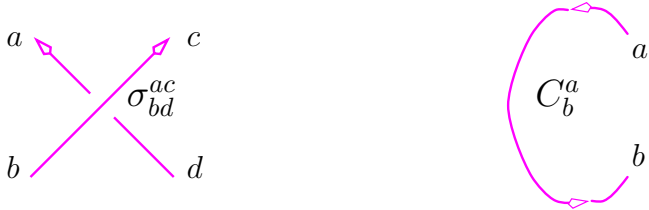}
    \caption{
      Diagram components corresponding to $\sigma$ and $C$.
    }
    \label{fig:SigmaC}
  \end{center}
\end{figure}

A \emph{state model} $(\sigma,C)$ for a link invariant consists of two
parameters: $\sigma$, an invertible rank $4$ tensor representing the
braid generator (i.e. a positive crossing), and $C$, an invertible rank
$2$ tensor representing a \emph{positive handle}, that is an
anticlockwise-oriented, vertical open arc, used to close a one string
of a braid. From these, we may immediately define the representations
corresponding to negative crossings (viz
$\overline{\sigma}\equiv\sigma^{-1}$) and negative handles (viz
$\overline{C}\equiv C^{-1}$).%
\footnote{%
  We frequently use the notation $\overline{X}$ to mean $X^{-1}$, in
  particular, writing $\overline{\sigma}\equiv \sigma^{-1}$ and
  $\overline{C}\equiv C^{-}$ allows us to omit superfluous ``$+$''
  signs, viz we write $C\equiv C^{+}$ for the positive handle.
}
Our current collection of arcs is shown in Figure \ref{fig:SigmaC}.

Let $\beta$ be a braid corresponding to a link $L\equiv\hat{\beta}$,
formed from the vertical closure of $\beta$.  The diagram components
corresponding to $\sigma$ (and $\overline{\sigma}$) are sufficient
to construct $\beta$, and those corresponding to $C$ (we don't need
$\overline{C}$) are then sufficient to construct $\hat{\beta}$ from
$\beta$.

When the pair $(\sigma,C)$ is chosen to satisfy the Reidemeister moves
(below, we write R1, R2 and R3), we may form a link invariant from the
contraction over the free indices of the tensors corresponding to the
diagram components.  It may happen that the algebraic structures
underlying the model mean that this invariant will be zero on closed
links (i.e.  $(0,0)$-tangles), however, we may still form an invariant
of $(1,1)$-tangles
\cite{AkutsuDeguchiOhtsuki:92,DeWitKauffmanLinks:99a}, by contracting
over all but one free index, and obtain an invariant which is not
necessarily trivial.  Our invariants $LG$ are based on typical
$U_q[gl(m|n)]$ representations, for which the appropriate supertrace is
zero, hence we define our invariants to be $(1,1)$-tangle invariants.


\subsection{State model parameters for $U_q[gl(m|n)]$ representations
            $\Lambda$}

Here, we integrate the materials of \S\ref{sec:quantumsuperalgebra} and
\S\ref{sec:Statemodelparameters}, allowing us to describe the
construction of state model parameters corresponding to arbitrary
$U_q[gl(m|n)]$ representations $\Lambda$.  Below, in \S\ref{sec:qhrho},
we perform some extra necessary calculations.  After that, in
\S\ref{sec:TheQuantumLinkInvariantsLGmn}, we specialise this material
to the case $\Lambda=(\dot{0}_m\,|\,\dot{\alpha}_n)$.

So, how do we construct state model parameters $(\sigma,C)$
corresponding to an invariant associated with an arbitrary
$U_q[gl(m|n)]$ representation $\pi$?

Firstly, the tensor product representation
$\check{R}\equiv(\pi\otimes\pi)\check{R}$ necessarily satisfies the
QYBE in the form (\ref{eq:QYBEinBraidForm}), and hence the braid
relation (\ref{eq:BraidRelation}). This means that abstract tensors
built from $\check{R}$ are invariant under R2 and R3, hence we may
construct representations of arbitrary braids from $\check{R}$.  Thus
$\sigma\triangleq\kappa_{\sigma}\check{R}$ (for any scalar constant
$\kappa_{\sigma}$) realises a representation of the braid generator.

A technical point distinguishes the quantum superalgebra situation
from that of the quantum algebra.  Quantum superalgebra R matrices are
in fact \emph{graded}, and actually satisfy a \emph{graded} QYBE.  It
is however, a simple matter to strip out this grading (i.e. apply an
automorphism \cite{DeWit:99d}),%
\footnote{%
  The $\check{R}$ supplied in \cite{DeWit:99d} are normalised such that
  $\lim_{q\to1}\check{R}$ is a (graded) permutation matrix. Scaling by
  $\kappa_{\sigma}$ does not change that.
}
yielding $\check{R}$ satisfying the usual, ungraded QYBE.

Secondly, to ensure that our invariant is an invariant of ambient
isotopy, we must select $C$ to ensure that abstract tensors built from
$\sigma$ and $C$ are also invariant under R1. To this end, we apply (a
grading-stripped version of) the following result
\cite[Lemma~2]{LinksGouldZhang:93} (see also \cite{LinksGould:92b}):
\begin{eqnarray*}
  ( I \otimes \mathrm{str} )
  [ ( I \otimes q^{2 h_\rho} ) \sigma ]
  =
  K I,
\end{eqnarray*}
where the Cartan element $q^{h_{\rho}}$ is defined in
\S\ref{sec:qhrho}, $\mathrm{str}$ is the supertrace, and $K$ is some
constant depending on the normalisations of $\sigma$ and
$q^{h_{\rho}}$.  Writing $S\equiv\pi(q^{2 h_\rho})$ for convenience;
for any scalar constant $\kappa_{C}$, setting $C\triangleq\kappa_{C}S$
allows us to represent positive handles.  It remains to choose
$(\kappa_S,\kappa_C)$ to satisfy R1.

\begin{figure}[htbp]
  \begin{center}
    \input{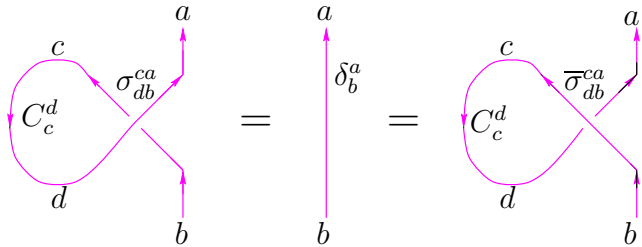}
    \caption{
      The first Reidemeister move.
    }
    \label{fig:Reidemeister1}
  \end{center}
\end{figure}

Thus, we demonstrate how to select $\kappa_{\sigma}$ and $\kappa_{C}$
such that the abstract tensor associated with removal of an isolated
loop is invariant under R1.  Figure \ref{fig:Reidemeister1} shows that
for $\sigma$ and $C$ to satisfy R1, they must satisfy (Einstein
summation convention):
\begin{equation}
  C^d_c
  \cdot
  \sigma^{c a}_{d b}
  =
  \delta^a_b
  =
  C^d_c
  \cdot
  \overline{\sigma}^{c a}_{d b},
  \label{eq:R1Constraints}
\end{equation}
where the definitions of $\kappa_{\sigma}$ and $\kappa_{C}$ yield:
$\overline{\sigma}=\kappa_{\sigma}^{-1}\check{R}^{-1}$,
and $\overline{C}=\kappa_{C}^{-1}S^{-1}$.

So, if we have established $S$ and $\check{R}$,
we may determine $\kappa_{\sigma}$ and $\kappa_{C}$ by solving
the following equations:
\begin{eqnarray*}
  \kappa_{C}
  \kappa_{\sigma}
  S^d_c
  \cdot
  (\check{R})^{c a}_{d b}
  =
  \kappa_{C}
  \kappa_{\sigma}^{-1}
  S^d_c
  \cdot
  (\check{R}^{-1})^{c a}_{d b}
  =
  \delta^a_b.
\end{eqnarray*}
Setting $a=b=1$ and using the fact that $S$ is diagonal, we thus have:
\begin{eqnarray*}
  \kappa_{\sigma}
  =
  X_1^{-\frac{1}{2}}
  X_2^{+\frac{1}{2}},
  \qquad
  \kappa_{C}
  =
  X_1^{-\frac{1}{2}}
  X_2^{-\frac{1}{2}},
\end{eqnarray*}
where $X_1 \triangleq S^c_c \cdot (\check{R})^{c 1}_{c 1}$ and
$X_2 \triangleq S^c_c \cdot (\check{R}^{-1})^{c 1}_{c 1}$.

Note that reflecting the diagrams of Figure \ref{fig:Reidemeister1}
about a vertical axis yields exactly the same constraints on
$\kappa_{\sigma}$ and $\kappa_{C}$.  To see this, the constraints
obtained by reflecting the diagrams in a vertical axis are:
\begin{equation}
  \overline{\sigma}^{a c}_{b d}
  \cdot
  \overline{C}^d_c
  =
  \delta^a_b
  =
  \sigma^{a c}_{b d}
  \cdot
  \overline{C}^d_c,
  \label{eq:ReflectedR1Constraints}
\end{equation}
however, we have:
$
  \overline{\sigma}^{a c}_{b d}
  =
  (\sigma^{c a}_{d b})|_{q\mapsto\overline{q}}
$
and $\overline{C}^d_c = (C^d_c)|_{q\mapsto\overline{q}}$.  Replacing
$q\mapsto\overline{q}$ in (\ref{eq:ReflectedR1Constraints}) and
applying these equivalences recovers (\ref{eq:R1Constraints}).
Similarly, reversing the orientations of the strings in Figure
\ref{fig:Reidemeister1} yields no new constraints.

What is significant in the above is that we have explicit formulae for
automatically scaling from $(\check{R},S)$ to $(\sigma,C)$, something
apparently absent in the literature. Variations on these formulae should
hold for a much wider class of representations and algebraic
structures. We write them up as a little lemma:

\begin{lemma}
  Let $\pi$ be a finite-dimensional highest weight $U_q[gl(m|n)]$
  representation, for which we have computed
  $\check{R}\equiv(\pi\otimes\pi)\check{R}$ and
  $S\equiv\pi(q^{2h_\rho})$. Then the state model parameters
  $(\sigma,C)$ for the corresponding link invariant of ambient isotopy 
  may be obtained from $(\check{R},S)$ by the scalings
  $\sigma=(X_1^{-1} X_2)^\frac{1}{2}\check{R}$ and
  $C=(X_1^{-1}X_2^{-1})^\frac{1}{2}S$, where
  $X_1\triangleq S^c_c\cdot(\check{R})^{c 1}_{c 1}$ and
  $X_2\triangleq S^c_c \cdot (\check{R}^{-1})^{c 1}_{c 1}$.
\end{lemma}


\subsubsection{Negative Handles, Caps and Cups}

Demanding that our model parameters satisfy R0 (ambient isotopy in the
plane) allows us to determine appropriate values for negative handles,
caps and cups (see \cite{DeWit:99a}). Although we \emph{can} evaluate
our invariants without these, we describe them here for completeness
and backwards compatibility.

Firstly, the negative handle $\overline{C}$ is simply
$C|_{q\mapsto\overline{q}}$.  Secondly, although there is some
flexibility in the choice of suitable caps $\Omega^{\pm}$ and cups
$\mho^{\pm}$, in fact it is natural to choose them to be the square
roots of the handles $C^{\pm}$:
\begin{equation}
  \Omega^{\pm}
  =
  \mho^{\pm}
  =
  (C^{\pm})^{\frac{1}{2}},
  \label{eq:DefnofCapsCups}
\end{equation}
taking the positive square root by convention.
Note that these choices further improve those of our previous work
\cite{DeWit:99a,DeWitKauffmanLinks:99a} by increasing the symmetries
between the diagram components.

\begin{figure}[htbp]
  \begin{center}
    \input{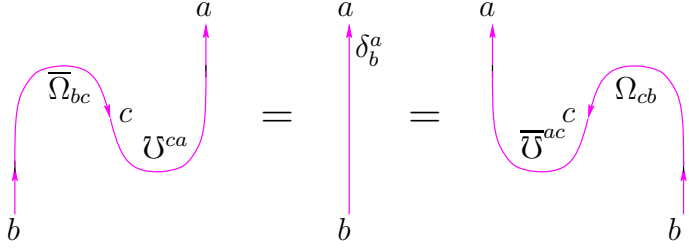}
    \caption{
      The zeroth Reidemeister move.
    }
    \label{fig:Reidemeister0}
  \end{center}
\end{figure}

Satisfaction of R0 is described in Figure \ref{fig:Reidemeister0}, that
is, we demand:
\begin{equation}
  \overline{\Omega}_{b c} \cdot \mho^{c a}
  =
  \delta^a_b
  =
  \overline{\mho}^{a c} \cdot \Omega_{c b}.
  \label{eq:R0Constraints}
\end{equation}

The definition (\ref{eq:DefnofCapsCups}) ensures that
(\ref{eq:R0Constraints}) is satisfied.  In fact, the LHS and RHS of
(\ref{eq:R0Constraints}) are actually equivalent, hence one is
redundant.  Again, reversing the orientations of the strings in Figure
\ref{fig:Reidemeister0} yields no new constraints.


\subsection{$q^{h_{\rho}}$ for $U_q[gl(m|n)]$}
\label{sec:qhrho}

Here, we determine the form that $q^{h_{\rho}}$ takes in
$U_q[gl(m|n)]$, in terms of Cartan generators.  Recall that for any
particular representation $\pi$, our state model requires (a
grading-stripped version of) $S=\pi(q^{2 h_{\rho}})$, and this may be
obtained by substitution of the appropriate matrix elements into the
expression for $q^{2 h_{\rho}}$.

Initially, we shall work with $gl(m|n)$. To this end, let $H$ be the
Cartan subalgebra of $gl(m|n)$, with dual the root space $H^{*}$. A
basis for $H^{*}$ is given by the \emph{fundamental weights}
$\{\varepsilon_i\}_{i=1}^{m+n}$, which are elementary unit vectors of
$m+n$ components, with $1$ in position $i$ and $0$ elsewhere.
On $H^*$, we have the following invariant bilinear form
$(\cdot,\cdot):H^{*} \times H^{*} \to \mathbb{C}$:
\begin{equation}
  (\varepsilon_i, \varepsilon_j)
  \triangleq
  (-)^{[i]} \delta_{i j},
  \label{eq:bilinearformonrootspace}
\end{equation}
and as $H$ and $H^{*}$ are dual, we of course have the form:
\begin{equation}
  {E^j}_j (\varepsilon_i)
  \triangleq
  \delta_{i j},
  \label{eq:actionofgensonrootspace}
\end{equation}
for $gl(m|n)$ Cartan generators ${E^j}_{j}$, $j=1,\dots,m+n$.

To the $gl(m|n)$ root $\varepsilon_i-\varepsilon_j$, there corresponds
a $gl(m|n)$ Chevalley generator ${E^i}_j$, and we assign a grading and
a sign to the roots in accordance with those of these generators.

In terms of these, $gl(m|n)$ has the following simple, positive roots:
\begin{equation}
  \alpha_i
  \triangleq
  \varepsilon_{i}
  -
  \varepsilon_{i+1},
  \qquad
  i = 1, \dots, m+n-1,
  \label{eq:Defnofroots}
\end{equation}
in the sense that these form a basis for $H^{*}$.  Apart from the
single odd root $\alpha_m$, the simple positive roots are all even.
(Of various choices for superalgebra root systems, this
\emph{distinguished root system} is unique in containing only one
odd root.)

Where $\Delta^{+}$ is the set of \emph{all} positive roots, and
$\gamma$ denotes the grading of the root $\gamma$, we define
$\rho$ as the graded half sum of all positive roots:
$
  \rho
  \triangleq
  \frac{1}{2}
  \sum_{\gamma\in\Delta^{+}}
    {(-)}^{[\gamma]}
    \gamma
$.
Explicitly, for $gl(m|n)$, we have
\cite[p6207]{DeliusGouldLinksZhang:95b}:
\begin{eqnarray*}
  \rho
  =
  {\textstyle \frac{1}{2}}
  \sum_{i=1}^{m}
    (m-n-2i+1)
    \varepsilon_{i}
  +
  {\textstyle \frac{1}{2}}
  \sum_{i=m+1}^{m+n}
    (3m+n-2i+1)
    \varepsilon_{i},
\end{eqnarray*}
although we will not actually require this form.

\vfill

\pagebreak

We are actually interested in $h_{\rho}\in gl(m|n)$, defined to satisfy:
\begin{equation}
  h_{\rho} (\alpha_i)
  \triangleq
  (\rho, \alpha_i),
  \qquad
  \forall \alpha_i,
  \label{eq:Defnofhrho}
\end{equation}
where we intend (\ref{eq:actionofgensonrootspace}) on the LHS and
(\ref{eq:bilinearformonrootspace}) on the RHS.
From the definition of $\rho$:
\begin{eqnarray}
  (\rho, \alpha_i)
  =
  {\textstyle \frac{1}{2}}
  (\alpha_i, \alpha_i)
  \stackrel{(\ref{eq:Defnofroots})}{=}
  {\textstyle \frac{1}{2}}
  (
    \varepsilon_{i} - \varepsilon_{i+1},
    \varepsilon_{i} - \varepsilon_{i+1}
  )
  \stackrel{(\ref{eq:bilinearformonrootspace})}{=}
  {\textstyle \frac{1}{2}}
  [
    (-)^{[i]} + (-)^{[i+1]}
  ].
  \label{eq:Dragon}
\end{eqnarray}
As $h_{\rho}$ is a Cartan element of $gl(m|n)$, we may
express it as a linear combination of Cartan generators
${E^i}_{i}$; viz for some undetermined scalar coefficients $\beta_i$,
we may set:
$
  h_{\rho}
  =
  \sum_{i=1}^{m+n}
    \beta_i
    {E^i}_i
$.
Substituting this into (\ref{eq:Defnofhrho}) yields:
\begin{eqnarray}
  h_{\rho} (\alpha_i)
  =
  \sum_{j=1}^{m+n}
    \beta_j
    {E^j}_j
    (\varepsilon_{i} - \varepsilon_{i+1})
  \stackrel{(\ref{eq:Defnofroots}, \ref{eq:actionofgensonrootspace})}{=}
  \beta_{i} - \beta_{i+1}.
  \label{eq:Tiger}
\end{eqnarray}
Substituting (\ref{eq:Dragon}) and (\ref{eq:Tiger}) into
(\ref{eq:Defnofhrho}), we have:
\begin{equation}
  \beta_{i} - \beta_{i+1}
  =
  \left\{
    \begin{array}{rl}
      + 1 & i = 1, \dots, m-1 \\
        0 & i = m \\
      - 1 & i = m+1, \dots, m+n-1.
    \end{array}
  \right.
  \label{eq:Ox}
\end{equation}
For symmetry, selecting $\beta_m = \theta$ and
substituting backwards and forwards yields:
\begin{eqnarray*}
  \beta_i
  =
  \theta
  +
  \left\{
    \begin{array}{ll}
      m-i     & i = 1, \dots, m \\
      i-(m+1) & i = m+1, \dots, m+n,
    \end{array}
  \right.
\end{eqnarray*}
therefore:
\begin{eqnarray*}
  h_{\rho}
  & = &
  \sum_{i=1}^{m+n}
    \beta_i
    {E^i}_i
  =
  \sum_{i=1}^{m}
    ( \theta + m-i) {E^i}_i
  +
  \sum_{i=m+1}^{m+n}
    ( \theta + i-(m+1)) {E^i}_i
  \\
  & = &
  \theta
  C_1
  +
  \sum_{i=1}^{m}
    (m-i) {E^i}_i
  +
  \sum_{i=m+1}^{m+n}
    (i-(m+1)) {E^i}_i,
\end{eqnarray*}
where $C_1\triangleq \sum_{i=1}^{m+n} {E^i}_i$ is the first-order
Casimir element of $gl(m|n)$. This shows us that $h_{\rho}$ is only
determined up to an additive constant.%
\footnote{%
  For $sl(m|n)$ and $sl(n)$, $h_{\rho}$ is actually unique.
  $C_1$ also satisfies $C_1(\alpha_i)=0, \forall \alpha_i$.
}

In passing from $gl(m|n)$ to $U_q[gl(m|n)]$, we pass from
$h_{\rho}$ to $q^{h_{\rho}}$, hence we have:
\begin{eqnarray*}
  q^{h_{\rho}}
  & = &
  q^{\theta C_1}
  \cdot
  q^{\sum_{i=1}^{m} (m-i) {E^i}_i}
  \cdot
  q^{\sum_{i=m+1}^{m+n} (i-(m+1)) {E^i}_i,}
  \\
  & = &
  {(q^{C_1})}^{\theta}
  \cdot
  \prod_{i=1}^{m} {\left( q^{{E^i}_i} \right)}^{m-i}
  \cdot
  \prod_{i=m+1}^{m+n} {\left( q^{{E^i}_i} \right)}^{i-(m+1)}
  \\
  & = &
  {(q^{C_1})}^{\theta}
  \cdot
  \prod_{i=1}^{m} K_i^{m-i}
  \cdot
  \prod_{i=m+1}^{m+n} K_i^{(m+1)-i},
\end{eqnarray*}
where we have reminded ourselves of the definition
$K_i \triangleq q^{(-)^{[i]} {E^i}_i}$.  Thus, $q^{h_{\rho}}$ is only
determined up to an arbitrary multiplicative constant.  Selecting
$\theta=0$, we declare the resulting product to be the standard
$q^{h_{\rho}}$. For arbitrary $m,n$, we have:
\begin{equation}
  q^{h_{\rho}}
  =
  K_1^{m-1}
  K_2^{m-2}
  \cdots
  K_{m-1}^{1}
  K_{m}^{0}
  \cdot
  K_{m+1}^{0}
  K_{m+2}^{-1}
  \cdots
  K_{m+n}^{-(n-1)},
  \label{eq:Mouse}
\end{equation}
where of course $K_i^0$ is the $U_q[gl(m|n)]$ identity element.

\pagebreak

For our state models we require $S=\pi(q^{2 h_{\rho}})$. To construct
$S$, it suffices to compute matrix elements for the $U_q[gl(m|n)]$
Cartan generators $K_i$, and insert (appropriate powers of) these into
(\ref{eq:Mouse}), finally stripping the grading from $S$.  In
\cite{DeWit:99d}, we described the automation of the construction of
$\check{R}$ corresponding to the $U_q[gl(m|1)]$ representations
$(\dot{0}_m\,|\,\alpha)$, for arbitrary $m$, and obtained explicit
$\check{R}$ for $m=1,2,3,4$. Explicit matrix elements for the $K_i$ are
obtained as a byproduct of that construction, facilitating the
evaluation of $S$.

What is particularly interesting about this work is that the entire
process, from the construction of the underlying representations
\cite{DeWit:99c,DeWit:99d}, to the scaling of the state model
parameters, to the final evaluations of the polynomials, has been
automated. This represents a step forward in computational power in
knot theory.


\section{The quantum link invariants $LG^{m,n}$}
\label{sec:TheQuantumLinkInvariantsLGmn}

Having described the construction of state models for arbitrary finite
dimensional highest weight $U_q[gl(m|n)]$ representations $\Lambda$, we
now restrict our attention to the case:
\begin{eqnarray*}
  \Lambda
  =
  (\dot{0}_m \; | \; \dot{\alpha}_n)
  \equiv
  (0, \dots, 0 \; | \; \alpha, \dots, \alpha),
\end{eqnarray*}
and the resulting invariants $LG^{m,n}$. Evaluation of $LG^{m,n}$ for
any particular link follows from that for $LG^{2}$, described in our
previous work \cite{DeWit:99a,DeWitKauffmanLinks:99a}.  Below, we make
a few comments on the properties of $LG^{m,n}$, before describing in
\S\ref{sec:Computationalissues} some computational issues and
evaluations for $LG^{3}$ and $LG^{4}$.


\subsection{Checking the QYBE and applying the
            Matveev $\Delta$--$\nabla$ test}

\begin{figure}[htbp]
  \begin{center}
    \input{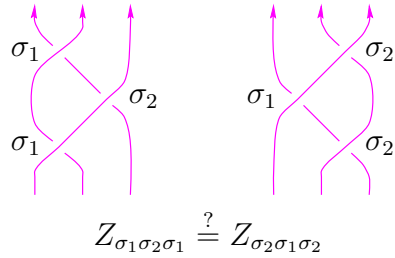}
    \caption{
      Checking that $\sigma$ satisfies the QYBE.
    }
    \label{fig:QYBECheck}
  \end{center}
\end{figure}

To be certain that we have made no errors in our computations, we check
that our braid generator $\sigma$ satisfies the (quantum) Yang--Baxter
equation. The code used to construct the tensors $Z_K$ is immediately
adaptable to such a test.  If $Z$ is the same for the braids
$\sigma_{1} \sigma_{2} \sigma_{1}$ and
$\sigma_{2} \sigma_{1} \sigma_{2}$,
then our braid generator satisfies the QYBE.
This is depicted in Figure \ref{fig:QYBECheck}.

\vspace{5mm}

\pagebreak

\begin{figure}[htbp]
  \begin{center}
    \input{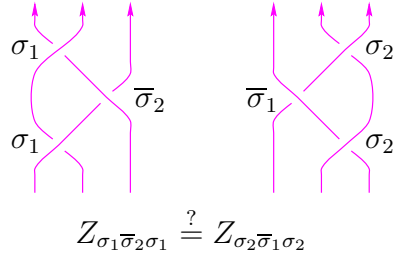}
    \caption{
      The Matveev $\Delta$--$\nabla$ test.
    }
    \label{fig:MatveevTest}
  \end{center}
\end{figure}

The same framework allows us to carry out a simple sufficiency check to
determine if a link invariant associated with some R matrix solution of
the QYBE will be trivial.%
\footnote{%
  This test is known to be a sufficient (but perhaps not a necessary)
  test of triviality -- it doesn't even guarantee the \emph{existence}
  of an invariant.
}
Matveev \cite{Matveev:87} (see also \cite{MurakamiNakanishi:89})
introduced a `delta unknotting operation' (which we call the Matveev
$\Delta$--$\nabla$ test), and proved that any knot can be transformed
to the unknot by using only this operation.  In our tensor language, if
$Z$ fails to distinguish $\sigma_{1}\overline{\sigma}_{2}\sigma_{1}$
and $\sigma_{2}\overline{\sigma}_{1}\sigma_{2}$, then the associated
invariant will be trivial, as a series of exchanges of crossings of
this form is always sufficient to convert any links to the unknot.
Matveev's test is depicted in Figure \ref{fig:MatveevTest}.

\vspace{5mm}

Both these tests have been satisfactorily carried out for our various
braid generators $\sigma$, viz each $\sigma$ satisfies the QYBE and the
invariant built from it is not necessarily trivial.


\subsection{Behaviour of $LG^{m,n}$ under inversion of $q$}

Let $K^{*}$ denote the reflection of a link $K$.  In
\cite{DeWitKauffmanLinks:99a}, we showed that
$LG^2_{K^*}=LG^2_{K}|_{q\to\overline{q}}$. This result immediately
carries over to $LG^{m,n}$, and means that if $LG^{m,n}_K$ is
palindromic in $q$ (i.e. invariant under the inversion
$q\to\overline{q}$), then $LG^{m,n}$ cannot distinguish the chirality
of $K$. Examples illustrating that $LG^{2}$ can distinguish the
chirality of all prime knots of up to $10$ crossings \cite{DeWit:99a}
demonstrate that $LG^{m,n}$ can indeed sometimes distinguish chirality,
although counterexamples are expected to exist.


\subsection{$LG^{m,n}$ doesn't detect mutation}

Theorem 5 of \cite{MortonCromwell:96} shows that quantum invariants
based on R matrices where the orthogonal decomposition of $V\otimes V$
contains no multiplicities will not distinguish mutants. The extension
of this result to quantum superalgebras is straightforward, and as our
invariants $LG^{m,n}$ are indeed based on representations of this type
\cite{DeliusGouldLinksZhang:95b}, they will not distinguish mutants.

\vfill

\pagebreak


\subsection{Behaviour of $LG^{m,n}$ under representation duality}

In \cite[Proposition~3.2]{DeWitKauffmanLinks:99a}, we showed that link
invariants derived from irreducible representations of quantum
(super)algebras are unable to detect knot inversion, as such invariants
are necessarily equivalent to invariants constructed from their dual
representations. Let us determine what this means for $LG^{m,n}$.

Let $V\equiv V_{\Lambda}$ be the module associated with
$\pi\equiv\pi_{\Lambda}$, viz $V$ has a highest weight vector $v_{+}$,
of weight $\Lambda$.  The corresponding lowest weight vector of $V$ is
obtained by the combined action of all the odd lowering operators:
$\prod {E^{m+j}}_{i}$ on $v_{+}$, where the product is over all
$i=1,\dots,m$ and $j=1,\dots,n$.  The action of ${E^{m+j}}_{i}$ on a
weight vector lowers its weight by $\varepsilon_{i}-\varepsilon_{m+j}$,
viz:
\begin{eqnarray*}
  (0, \dots, 0, 1, 0, \dots, 0 \; | \; 0, \dots, 0, -1, 0, \dots, 0).
\end{eqnarray*}
The resulting lowest weight vector of $V$ thus has weight
$\overline{\Lambda}$:
\begin{eqnarray*}
  \overline{\Lambda}
  & = &
  \Lambda
  -
  {\textstyle \sum_{i j}}
    (\varepsilon_{i}-\varepsilon_{m+j})
  =
  \Lambda
  -
  n
  {\textstyle \sum_{i=1}^{m}}
    \varepsilon_{i}
  -
  m
  {\textstyle \sum_{j=1}^{n}}
    \varepsilon_{m+j}
  \\
  & = &
  (-n, \dots, -n \; | \; \alpha+m, \dots, \alpha+m).
\end{eqnarray*}
The dual of $V$ is labeled $V^*$, and naturally has highest weight
$-\overline{\Lambda}$:
\begin{eqnarray*}
  - \overline{\Lambda}
  =
  (n, \dots, n \; | \; -\alpha-m, \dots, -\alpha-m),
\end{eqnarray*}
but $V^{*}$ is equivalent to the module of highest weight $\Lambda^{*}$:
\begin{eqnarray*}
  \Lambda^{*}
  =
  (0, \dots, 0 \; | \; -\alpha-(m-n), \dots, -\alpha-(m-n)),
\end{eqnarray*}
hence we may regard the representations $\Lambda$ and $\Lambda^*$ as
duals.  Thus, at least up to a scalar multiple, we expect $LG^{m,n}$ to
be invariant under the transformation $\alpha \mapsto -\alpha-(m-n)$,
equivalently in the more symmetric form:
$\alpha +\frac{m-n}{2}\mapsto -\alpha-\frac{m-n}{2}$, viz:
\begin{eqnarray*}
  LG^{m,n}(q, q^{ \alpha+\frac{m-n}{2}})
  =
  LG^{m,n}(q, q^{-\alpha-\frac{m-n}{2}}).
\end{eqnarray*}
Thus, if we define:
\begin{equation}
  p
  \triangleq
  q^{ \alpha+\frac{m-n}{2}},
  \label{eq:Defnofp}
\end{equation}
we have, again, up to a scalar multiple, the symmetry:
\begin{equation}
  LG^{m,n}(q, p)
  =
  LG^{m,n}(q, \overline{p}).
  \label{eq:LGsymmetryunderinversionofp}
\end{equation}
Experiments show that the scalar multiple is always $\pm 1$,
and, for knots, always $1$.
Where $\overline{K}$ is the inverse of a knot $K$,
inspection of diagram components shows that
$
  LG^{m,n}_{\overline{K}}(q, p)
  =
  LG^{m,n}_{K}(q, \overline{p})
$,
hence (\ref{eq:LGsymmetryunderinversionofp}) shows that $LG^{m,n}$ is
unable to detect the inversion of knots.

Experimentally (setting $n=1$), we find that the only time the ``$-$''
sign actually appears is for odd $m$ and links of $2$ components.  That
this is true for case $m=1$ (which is in fact the Alexander--Conway
polynomial) is well-known \cite{Lickorish:97}.  These results are
exemplified in our previous work
\cite{DeWit:99a,DeWitKauffmanLinks:99a} for $LG^2$.

Lastly, we often wish to eliminate $\alpha$ from expressions of the form
$q^{x\alpha+y}$, to express them in terms of $p$ and $q$ alone. Using
(\ref{eq:Defnofp}), we have:
\begin{eqnarray*}
  q^{x\alpha+y}
  =
  p^x q^{y - x (\frac{m-n}{2})}.
\end{eqnarray*}


\subsection{$LG^{m,n}$ of split links}

Recall that we define $LG^{m,n}_K$ as a $(1,1)$ tangle invariant,
obtained for a link $K$ as the first component of the diagonal tensor
(scalar multiple of the identity) $T_K$. We do this as the closed form
(i.e. the $(0,0)$ tangle form) always evaluates to zero (cf. the ADO
invariant \cite{AkutsuDeguchiOhtsuki:92}).

To see this, begin by observing that the value of our state model on
$0_1$ (i.e. the unknot, an isolated loop) as a $(0,0)$ tangle is zero,
as $\sum_a C^a_a=0$.  This follows from the fact that for the
$U_q[gl(m|n)]$ superalgebras, the $q$-superdimension of typical
representations (defined by $\mathrm{str}[\pi(q^{2 h_{\rho}})]$) is
always identically zero \cite{LinksGould:96b}.  As $S$ is a
grading-stripped version of the exponential of the Cartan element
$\pi(q^{2 h_{\rho}})$, we necessarily have $\mathrm{tr}(S)=0$, hence
$\mathrm{tr}(C)=0$.  Multiplying these results by the scalar in $T_K$
yields the result.

Now, let $K=K_1 \sqcup K_2$ be the split (i.e. disconnected, separated)
union of links $K_1$ and $K_2$, and say that we are trying to evaluate
the $(1,1)$ tangle form using a string of $K_1$.  The construction of
$LG^{m,n}_K$ means that at some point of contracting $Z_K$ to $T_K$, we
close the final string of $K_2$, and at this stage our tensor becomes
zero throughout, thus $T_K$ is zero.  Thus, as disconnected
multicomponent links represented by $(1,1)$ tangles
necessarily include a closed component, we
have proven:
\begin{theorem}
  $LG^{m,n}_K=0$ for disconnected multicomponent links $K$.
\end{theorem}


\section{Computational issues in evaluating $LG^{m,n}$}
\label{sec:Computationalissues}

\subsection{Various sets of computational variables}

The representation of the braid generator $\sigma$ obtained from
the representation theory \cite{DeWit:99c,DeWit:99d} contains
algebraic expressions in variables $q$ and $\alpha$, including
many $q$ brackets. This form is readable to
human eyes, but can be improved upon for machine consumption.
We shall call $\{q,\alpha\}$ the \textbf{rep(resentation)} variables.

From (\ref{eq:Defnofp}), we see that our link invariants are naturally
expressed in terms of $q$ and $p\triangleq q^{\alpha+\frac{m-n}{2}}$;
so we initially make this change of variables in the internal
representation of the braid generator and the positive handle. This
action replaces all the $q$ brackets, which contained $\alpha$.  The
resulting braid generator contains rational expressions in variables
$q^{\frac{1}{2}}$ and $p$.  To simplify the vulgar fractions within the
exponents, we define a new variable to be used internally:
$Q\triangleq q^{\frac{1}{2}}$.  In some sense, the resulting braid
generator is now optimally literate, and we use this form to accrete
tensors to build $Z_K$, and also to check the QYBE and the Matveev
$\Delta$--$\nabla$ test.  We shall call $\{Q,p\}$ the
\textbf{int(ernal)} variables, and to convert from \textbf{rep} to
\textbf{int} variables, we shall invoke \emph{in order} the following
rules:
\begin{eqnarray*}
  \left\{
    q^{x\alpha+y}
    \mapsto p^x q^{y - x
    (\frac{m-n}{2})},
    \qquad
    q \mapsto Q^2
  \right\},
\end{eqnarray*}
where $x,y\in\mathbb{Z}$.  We occasionally have an interest in the
inverse transformation to convert from \textbf{int} to \textbf{rep}
variables, and for this we shall invoke \emph{in order} the following
rules:
\begin{eqnarray*}
  \hspace{-3pt}
  \left\{
    Q \mapsto q^\frac{1}{2},
    \qquad
    p \mapsto
    q^{\alpha+\frac{m-n}{2}},
    \qquad
    (q^{x\alpha+y} - \overline{q}^{x\alpha+y})
    \mapsto
    (q - \overline{q}) [x\alpha+y]_q
  \right\},
\end{eqnarray*}
where handling the last of these rules typically requires some care.

Sometimes, we must invert the \textbf{int} variables, for example in
computing the inverse braid generator $\overline{\sigma}$.  We have the
rules:
\begin{eqnarray*}
  \left\{
    Q\mapsto \overline{Q},
    \qquad
    p \mapsto \overline{p}
  \right\}.
\end{eqnarray*}

Finally, extracting the first component of $T_K$ thus yields an
expression \textbf{int} variables. We must then expand $Q\mapsto q^2$.
Furthermore, we discover that $LG^{m,n}_K$ is actually an invariant in
$p^2$ not just $p$, so we define $P=p^2$ to reduce things a little.  We
shall call $\{q,P\}$ the \textbf{L(ink) I(nvariant)} variables, and to
convert from \textbf{int} to \textbf{LI} variables, we shall invoke the
following rules:
\begin{eqnarray*}
  \left\{
    Q \mapsto q^{\frac{1}{2}},
    \qquad
    p \mapsto P^{\frac{1}{2}}
  \right\}.
\end{eqnarray*}

Parameters used for the state models for $LG^m$ for $m=1,2,3,4$ are
presented in Appendix \ref{app:StateModelParameters}.


\subsection{Explicit construction of $S$}

From (\ref{eq:Mouse}), we have for $U_q[gl(m|n)]$ that
$S\equiv \pi(q^{2 h_{\rho}})$ is:
\begin{eqnarray}
  S
  & = &
  \pi(K_1)^{2(m-1)}
  \pi(K_2)^{2(m-2)}
  \cdots
  \pi(K_{m-1})^{2}
  \cdot
  \nonumber
  \\
  & &
  \qquad
  \pi(K_{m+2})^{-2}
  \pi(K_{m+3})^{-4}
  \cdots
  \pi(K_{m+n})^{-2(n-1)}.
  \label{eq:Rabbit}
\end{eqnarray}
Setting $n=1$ in (\ref{eq:Rabbit}), we have:
\begin{eqnarray*}
  S
  =
  \pi_{\Lambda} (K_1)^{2(m-1)}
  \pi_{\Lambda} (K_2)^{2(m-2)}
  \cdots
  \pi_{\Lambda} (K_{m-1})^{2}.
\end{eqnarray*}
To illustrate, for the $U_q[gl(2|1)]$ case, we have $h_{\rho}={E^1}_1$,
hence $q^{h_{\rho}}=K_1$, so $S=\pi(K_1)^2$.  This contrasts
with the choice of $\theta=-1$ made in \cite{DeWit:98}, which yields
$h_{\rho}=-{E^2}_2-{E^3}_3$, viz
$q^{h_{\rho}} = q^{-{E^2}_2} q^{-{E^3}_3} = K_2^{-1} K_3$, so
$S=\pi(K_2)^{-2} \pi(K_3)^2$.


\subsection{Illustrative examples of $LG^{m}$}

At present, we are able to compute state model parameters for
$LG^{m,1}\equiv LG^{m}$ only, as we have not yet computed
$\check{R}$ or matrix elements for the $K_i$ for cases $n\neq 1$.
For the cases $m=1,2,3,4$, we are able to make the following comments.

\begin{itemize}
\item
  $LG^1$ is the Alexander--Conway polynomial in variable
  $P\equiv q^{2\alpha}$.  This is a well-known result, cf.
  \cite{LiaoSong:91,RozanskySaleur:92}.

\item
  Evaluations for $LG^2$ for all prime knots of up to $10$ crossings
  have been reported in \cite{DeWit:99a}. In that paper, we claimed
  that $LG^m$ for $m>2$ was essentially incomputable due to vast memory
  requirements of the tensors $Z_K$; but we have since made some
  headway in this by adapting our code to recognise the sparsity of
  these tensors; doing the symbolic equivalent of what is called
  ``sparse matrix multiplication'' in numerical linear algebra.  This
  change comes at a cost of more lines of interpreted code, but is
  still an improvement in algorithmic efficiency.  It also results in
  an increase in the speed of computation for $LG^2$, and facilitates
  its evaluation from braid presentations of $6$ strings, something not
  previously feasible.

\item
  Evaluations for $LG^3$ and $LG^4$ for various links are
  presented in Appendix \ref{app:LG34}.  Those lists are quite brief,
  and only include some links of braid index at most $3$.  Our current
  computational method requires too much memory for us to extend our
  tables of polynomials any further.

\end{itemize}

Of some interest is the rate of growth in exponent of the polynomials
with $m$ for a particular link.  For example, we have the following
results for the trefoil knot $3_1$ and the figure eight knot $4_1$:

\begin{eqnarray*}
  \begin{array}{rcl}
    LG^1_{3_{1}}
    & = &
    -
    \hspace{44pt}
    (1)
    \\
    & &
    +
    (\overline{P}^{1} + P^{1}) (+ 1)
    \\
    \\
    LG^2_{3_{1}}
    & = &
    \hspace{52pt}
    (1 + 2 q^{2})
    \\
    & &
    -
    (\overline{P}^{1} + P^{1}) (q + q^{3})
    \\
    & &
    +
    (\overline{P}^{2} + P^{2}) (q^{2})
    \\
    \\
    LG^3_{3_{1}}
    & = &
    -
    \hspace{44pt}
    (q^{2} + 2 q^{4}  + 3 q^{6}  + q^{8})
    \\
    & &
    +
    (\overline{P}^{1} + P^{1}) (q^{2}  + 2 q^{4}  + 2 q^{6}  + q^{8})
    \\
    & &
    -
    (\overline{P}^{2} + P^{2}) (q^{4}  + q^{6}  + q^{8})
    \\
    & &
    +
    (\overline{P}^{3} + P^{3}) (q^{6})
    \\
    \\
    LG^4_{3_{1}}
    & = &
    \hspace{52pt}
    (q^{4} + 2 q^{6} + 4 q^{8} + 4 q^{10} + 5 q^{12} + 2 q^{14} +q^{16})
    \\
    & &
    -
    (\overline{P}^{1} + P^{1})
    (q^{5} + 2 q^{7} + 4 q^{9} + 4 q^{11} + 3 q^{13} + 2 q^{15})
    \\
    & &
    +
    (\overline{P}^{2} + P^{2})
    (  q^{6} + 2 q^{8} + 2 q^{10} + 3 q^{12} + q^{14} + q^{16})
    \\
    & &
    -
    (\overline{P}^{3} + P^{3})
    (q^{9} + q^{11} + q^{13} + q^{15})
    \\
    & &
    +
    (\overline{P}^{4} + P^{4})
    (  q^{12})
  \end{array}
\end{eqnarray*}

\begin{eqnarray*}
  \begin{array}{rcl}
    LG^1_{4_{1}}
    & = &
    \hspace{52pt}
    (3)
    \\
    & &
    -
    (\overline{P}^{1} + P^{1}) (1)
    \\
    \\
    LG^2_{4_1}
    & = &
    \hspace{52pt}
    (2 \overline{q}^{2} + 7 + 2 q^{2})
    \\
    & &
    -
    (\overline{P}^{1} + P^{1}) (3 \overline{q} + 3 q)
    \\
    & &
    +
    (\overline{P}^{2} + P^{2})
    (1)
    \\
    \\
    LG^3_{4_1}
    & = &
    \hspace{52pt}
    (5 \overline{q}^{4} + 9 \overline{q}^{2} + 17 + 9 q^{2} + 5 q^{4})
    \\
    & &
    -
    (\overline{P}^{1} + P^{1})
    (2 \overline{q}^{4} + 8 \overline{q}^{2} + 10 + 8 q^{2} + 2 q^{4})
    \\
    & &
    +
    (\overline{P}^{2} + P^{2}) (3 \overline{q}^{2} + 3 + 3 q^{2})
    \\
    & &
    -
    (\overline{P}^{3} + P^{3}) (1)
  \end{array}
\end{eqnarray*}

\pagebreak


\section{Further work}

The current work is part of a larger program to automate the
construction of more general quantum link invariants.
A few comments on the direction of this program are in order.

\begin{itemize}
\item
  In this paper, the limits of our method of evaluation have been
  reached,%
  \footnote{%
    The material has also been applied to the evaluation of `$N$-Jones'
    polynomials $V^N$. These are the quantum link invariants associated
    with the $N$ dimensional representations of $U_q[sl(2)]$.  In the
    language of \cite{KirbyMelvin:91}, they are monochromatic versions
    of coloured Jones polynomials of order $N$.  The limit to
    computation for these invariants for prime knots of up to $10$
    crossings is around $N=4$, although we can calculate
    $V^{13}_{3_1}$.
  }
  and a more efficient method of evaluation is required.  A promising
  candidate involves chasing through braids one crossing at a time,
  accumulating only an $N \times N$ matrix (where $N=\mathrm{dim(V)}$)
  of polynomials at each step.  That method requires foreknowledge of
  the decomposition of $\check{R}$ into the canonical form
  $\check{R}=\sum_ia_i\otimes b_i$, and this is already available for
  $U_q[sl(2)]$ and $U_q[gl(1|1)]$.  It is applicable to links of any
  number of crossings and components, and is really only limited by
  $N$, although much less strongly than our current method.  In
  particular, it is not dependent on the string index of braid
  presentations.

\item
  Moreover, the construction of more general quantum link invariants
  requires a more general approach to construction of underlying R
  matrices. The current method \cite{DeWit:99c,DeWit:99d} exploits
  explicit knowledge of the decomposition of the tensor product of the
  underlying module, but this is \emph{not} generally known.
  Alternatively, it is also possible to construct explicit R matrices
  from knowledge of the universal (i.e. representation-independent) R
  matrix and the matrix elements of the underlying representation.  As
  we have to hand details of the universal R matrices for
  arbitrary quantum (super)algebras \cite{KhoroshkinTolstoy:91} (albeit
  in a somewhat abstract form), and some knowledge of a process to
  construct the matrix elements, it is eminently possible to construct
  many more R matrices.

\item
  Lastly, we are limited by our use of braids, for which we have
  systematic tables only for the first $249$ prime knots of up to $10$
  crossings. As of 1998, Dowker codes for all the $1,701,936$ prime
  knots of up to $16$ crossings have been enumerated
  \cite{HosteThistlethwaiteWeeks:98}, and our not being able to access
  them is a sad thing. As we don't have the implementation of an
  algorithm that allows us to map these Dowker codes to braids, it
  is attractive to try to adapt new material to accept Dowker codes as
  input. The converse to this is that our new invariants $LG^{m,n}$
  are well suited to extending those tables, as they distinguish many
  more knots than other polynomial invariants.

\end{itemize}



\enlargethispage{\baselineskip}

\section*{Acknowledgements}

I am grateful to Jon Links of The University of Queensland, Australia,
for advice on how the deduce the natural choice of variables to
describe $LG^{m,n}$, and other worthwhile discussions.

My research at Kyoto University is funded by a Postdoctoral
Fellowship for Foreign Researchers (\# P99703), provided by the Japan
Society for the Promotion of Science.
D\={o}mo arigat\={o} gozaimashita!

\vfill

\pagebreak


\appendix

\section{State model parameters}
\label{app:StateModelParameters}

Below, we list state model parameters for $LG^{m}$, for $m=1,2,3,4$. To
improve literacy, we have written $[X]$ for $[X]_q$, $\overline{X}$ for
$X^{-1}$, for various $X$, and $\Delta = q - \overline{q}$.  Horizontal
lines divide tensor components into symmetry classes.

\subsection*{Parameters for $LG^{1}$}

The braid generator $\sigma$ has $5$ nonzero components:
\begin{eqnarray*}
  &&
  \hspace{-45pt}
  q^{- \alpha}
  \left\{
    \begin{array}{@{\hspace{0mm}}c@{\hspace{0mm}}}
      e^{1 1}_{1 1}
    \end{array}
  \right\},
  \qquad
  -
  q^{\alpha}
  \left\{
    \begin{array}{@{\hspace{0mm}}c@{\hspace{0mm}}}
      e^{2 2}_{2 2}
    \end{array}
  \right\},
  \qquad
  -
  \Delta
  [\alpha]
  \left\{
    \begin{array}{@{\hspace{0mm}}c@{\hspace{0mm}}}
        e^{2 1}_{2 1}
    \end{array}
  \right\},
  \qquad
  1
  \left\{
    \begin{array}{@{\hspace{0mm}}c@{\hspace{0mm}}}
        e^{1 2}_{2 1}
        \\
        e^{2 1}_{1 2}
    \end{array}
  \right\},
\end{eqnarray*}
and the left handle $C$ has $2$ diagonal components:
\begin{eqnarray*}
  C
  =
  q^{- \alpha}
  \left\{
    \begin{array}{l}
      + e^{1}_{1}
      \\
      - e^{2}_{2}
    \end{array}
  \right\},
\end{eqnarray*}
using the scaling factors:
\begin{eqnarray*}
  \kappa_{\sigma}
  =
  q^{- \alpha},
  \qquad \qquad
  \kappa_{C}
  =
  q^{- \alpha}.
\end{eqnarray*}

\subsection*{Parameters for $LG^{2}$}

The braid generator $\sigma$ has $26$ nonzero components:
\begin{eqnarray*}
  &&
  \hspace{-45pt}
  q^{-2 \alpha}
  \left\{
    \begin{array}{@{\hspace{0mm}}c@{\hspace{0mm}}}
      e^{1 1}_{1 1}
    \end{array}
  \right\},
  \qquad
  - 1
  \left\{
    \begin{array}{@{\hspace{0mm}}c@{\hspace{0mm}}}
      e^{2 2}_{2 2},
      e^{3 3}_{3 3}
    \end{array}
  \right\},
  \qquad
  q^{2 \alpha + 2}
  \left\{
    \begin{array}{@{\hspace{0mm}}c@{\hspace{0mm}}}
      e^{4 4}_{4 4}
    \end{array}
  \right\},
  \\[1mm]
  \hline
  \\[1mm]
  &&
  \hspace{-45pt}
  -
  \Delta
  q^{- \alpha}
  [\alpha]
  \left\{
    \begin{array}{@{\hspace{0mm}}c@{\hspace{0mm}}}
      e^{2 1}_{2 1},
      e^{3 1}_{3 1}
    \end{array}
  \right\},
  \qquad
  \Delta^2
  q
  [\alpha]
  [\alpha + 1]
  \left\{
    \begin{array}{@{\hspace{0mm}}c@{\hspace{0mm}}}
      e^{4 1}_{4 1}
    \end{array}
  \right\},
  \\
  &&
  \hspace{-45pt}
  \Delta
  q^{\alpha + 1}
  [\alpha + 1]
  \left\{
    \begin{array}{@{\hspace{0mm}}c@{\hspace{0mm}}}
      e^{4 2}_{4 2},
      e^{4 3}_{4 3}
    \end{array}
  \right\},
  \qquad
  \Delta
  q
  \left\{
    \begin{array}{@{\hspace{0mm}}c@{\hspace{0mm}}}
      e^{3 2}_{3 2}
    \end{array}
  \right\},
  \\[1mm]
  \hline
  \\[1mm]
  &&
  \hspace{-45pt}
  1
  \left\{
    \begin{array}{@{\hspace{0mm}}c@{\hspace{0mm}}}
      e^{1 4}_{4 1}
      \\
      e^{4 1}_{1 4}
    \end{array}
  \right\},
  \qquad
  -
  q
  \left\{
    \begin{array}{@{\hspace{0mm}}c@{\hspace{0mm}}}
      e^{2 3}_{3 2}
      \\
      e^{3 2}_{2 3}
    \end{array}
  \right\},
  \qquad
  q^{-\alpha}
  \left\{
    \begin{array}{@{\hspace{0mm}}c@{\hspace{0mm}}}
      e^{1 2}_{2 1},
      e^{1 3}_{3 1}
      \\
      e^{2 1}_{1 2},
      e^{3 1}_{1 3}
    \end{array}
  \right\},
  \qquad
  q^{\alpha + 1}
  \left\{
    \begin{array}{@{\hspace{0mm}}c@{\hspace{0mm}}}
      e^{2 4}_{4 2},
      e^{3 4}_{4 3}
      \\
      e^{4 2}_{2 4},
      e^{4 3}_{3 4}
    \end{array}
  \right\},
  \\[1mm]
  \hline
  \\[1mm]
  &&
  \hspace{-45pt}
  \Delta
  q
  [\alpha]^{\frac{1}{2}}
  [\alpha + 1]^{\frac{1}{2}}
  \left\{
    -
    \overline{q}^{\frac{1}{2}}
    \left\{
    \begin{array}{@{\hspace{0mm}}c@{\hspace{0mm}}}
      e^{2 3}_{4 1}
      \\
      e^{4 1}_{2 3}
    \end{array}
    \right\},
    +
    q^{\frac{1}{2}}
    \left\{
    \begin{array}{@{\hspace{0mm}}c@{\hspace{0mm}}}
      e^{3 2}_{4 1}
      \\
      e^{4 1}_{3 2}
    \end{array}
    \right\}
  \right\},
\end{eqnarray*}
and the left handle $C$ has $4$ diagonal components:
\begin{eqnarray*}
   C
   =
   q^{- 2 \alpha - 1}
   \left\{
     \begin{array}{l}
       +
       q
       \left\{
                      e^{1}_{1}
       \right\}
       \\
       -
       \left\{
         q            e^{2}_{2}, \; \;
         \overline{q} e^{3}_{3}
       \right\}
       \\
       +
       \overline{q}
       \left\{
                      e^{4}_{4}
       \right\}
     \end{array}
   \right\},
\end{eqnarray*}
using the scaling factors:
\begin{eqnarray*}
  \kappa_{\sigma}
  =
  q^{- 2 \alpha},
  \qquad \qquad
  \kappa_{C}
  =
  q^{- 2 \alpha}.
\end{eqnarray*}

\pagebreak

\subsection*{Parameters for $LG^{3}$}

The braid generator $\sigma$ has $139$ nonzero components:

\begin{eqnarray*}
  &&
  \hspace{-45pt}
  q^{- 3 \alpha}
  \left\{
    \begin{array}{@{\hspace{0mm}}c@{\hspace{0mm}}}
      e^{1 1}_{1 1}
    \end{array}
  \right\},
  \quad \hspace{-4pt}
  -
  q^{- \alpha}
  \left\{
    \begin{array}{@{\hspace{0mm}}c@{\hspace{0mm}}}
      e^{2 2}_{2 2},
      e^{3 3}_{3 3},
      e^{4 4}_{4 4}
    \end{array}
  \right\},
  \quad \hspace{-4pt}
  q^{\alpha + 2}
  \left\{
    \begin{array}{@{\hspace{0mm}}c@{\hspace{0mm}}}
      e^{5 5}_{5 5},
      e^{6 6}_{6 6},
      e^{7 7}_{7 7}
    \end{array}
  \right\},
  \quad \hspace{-4pt}
  -
  q^{3 \alpha + 6}
  \left\{
    \begin{array}{@{\hspace{0mm}}c@{\hspace{0mm}}}
      e^{8 8}_{8 8}
    \end{array}
  \right\},
  \\[1mm]
  \hline
  \\[1mm]
  &&
  \hspace{-45pt}
  -
  \Delta
  q^{- 2 \alpha}
  [\alpha]
  \left\{
    \begin{array}{@{\hspace{0mm}}c@{\hspace{0mm}}}
      e^{2 1}_{2 1},
      e^{3 1}_{3 1},
      e^{4 1}_{4 1}
    \end{array}
  \right\},
  \qquad
  \Delta^2
  q^{- \alpha + 1}
  [\alpha]
  [\alpha + 1]
  \left\{
    \begin{array}{@{\hspace{0mm}}c@{\hspace{0mm}}}
      e^{5 1}_{5 1},
      e^{6 1}_{6 1},
      e^{7 1}_{7 1}
    \end{array}
  \right\},
  \\
  &&
  \hspace{-45pt}
  -
  \Delta
  q^{2 \alpha + 4}
  [\alpha + 2]
  \left\{
    \begin{array}{@{\hspace{0mm}}c@{\hspace{0mm}}}
      e^{8 7}_{8 7},
      e^{8 6}_{8 6},
      e^{8 5}_{8 5}
    \end{array}
  \right\},
  \qquad
  -
  \Delta^2
  q^{\alpha + 3}
  [\alpha + 1]
  [\alpha + 2]
  \left\{
    \begin{array}{@{\hspace{0mm}}c@{\hspace{0mm}}}
       e^{8 4}_{8 4},
       e^{8 3}_{8 3},
       e^{8 2}_{8 2}
    \end{array}
  \right\},
  \\
  &&
  \hspace{-45pt}
  \Delta
  q
  [\alpha + 1]
  \left\{
    \begin{array}{@{\hspace{0mm}}c@{\hspace{0mm}}}
      e^{5 2}_{5 2},
      e^{5 3}_{5 3},
      e^{6 2}_{6 2},
      e^{6 4}_{6 4},
      e^{7 3}_{7 3},
      e^{7 4}_{7 4}
    \end{array}
  \right\},
  \\
  &&
  \hspace{-45pt}
  \Delta
  q^{- \alpha + 1}
  \left\{
    \begin{array}{@{\hspace{0mm}}c@{\hspace{0mm}}}
      e^{3 2}_{3 2},
      e^{4 2}_{4 2},
      e^{4 3}_{4 3}
    \end{array}
  \right\},
  \qquad
  -
  \Delta
  q^{\alpha + 3}
  \left\{
    \begin{array}{@{\hspace{0mm}}c@{\hspace{0mm}}}
      e^{6 5}_{6 5},
      e^{7 5}_{7 5},
      e^{7 6}_{7 6}
    \end{array}
  \right\},
  \\
  &&
  \hspace{-45pt}
  \Delta
  q^{2}
  [\alpha + 1]
  \left\{
    \begin{array}{@{\hspace{0mm}}c@{\hspace{0mm}}}
      - \Delta               e^{6 3}_{6 3},
      q (\overline{q}^2 - q^2) e^{7 2}_{7 2}
    \end{array}
  \right\},
  \qquad
  -
  \Delta^3
  q^{3}
  [\alpha]
  [\alpha + 1]
  [\alpha + 2]
  \left\{
    \begin{array}{@{\hspace{0mm}}c@{\hspace{0mm}}}
      e^{8 1}_{8 1}
    \end{array}
  \right\},
  \\[1mm]
  \hline
  \\[1mm]
  &&
  \hspace{-45pt}
  q^{- 2 \alpha}
  \left\{
    \begin{array}{@{\hspace{0mm}}c@{\hspace{0mm}}}
      e^{1 2}_{2 1},
      e^{1 3}_{3 1},
      e^{1 4}_{4 1}
      \\
      e^{2 1}_{1 2},
      e^{3 1}_{1 3},
      e^{4 1}_{1 4}
    \end{array}
  \right\},
  \qquad
  q^{- \alpha}
  \left\{
    \begin{array}{@{\hspace{0mm}}c@{\hspace{0mm}}}
      e^{1 5}_{5 1},
      e^{1 6}_{6 1},
      e^{1 7}_{7 1}
      \\
      e^{5 1}_{1 5},
      e^{6 1}_{1 6},
      e^{7 1}_{1 7}
    \end{array}
  \right\},
  \\
  &&
  \hspace{-45pt}
  q^{2 \alpha + 4}
  \left\{
    \begin{array}{@{\hspace{0mm}}c@{\hspace{0mm}}}
      e^{8 5}_{5 8},
      e^{8 6}_{6 8},
      e^{8 7}_{7 8}
      \\
      e^{5 8}_{8 5},
      e^{6 8}_{8 6},
      e^{7 8}_{8 7}
    \end{array}
  \right\},
  \qquad
  -
  q^{\alpha + 2}
  \left\{
    \begin{array}{@{\hspace{0mm}}c@{\hspace{0mm}}}
      e^{8 2}_{2 8},
      e^{8 3}_{3 8},
      e^{8 4}_{4 8}
      \\
      e^{2 8}_{8 2},
      e^{3 8}_{8 3},
      e^{4 8}_{8 4}
    \end{array}
  \right\},
  \\
  &&
  \hspace{-45pt}
  -
  q^{- \alpha + 1}
  \left\{
    \begin{array}{@{\hspace{0mm}}c@{\hspace{0mm}}}
      e^{2 3}_{3 2},
      e^{2 4}_{4 2},
      e^{3 4}_{4 3}
      \\
      e^{3 2}_{2 3},
      e^{4 2}_{2 4},
      e^{4 3}_{3 4}
    \end{array}
  \right\},
  \qquad
  q^{\alpha + 3}
  \left\{
    \begin{array}{@{\hspace{0mm}}c@{\hspace{0mm}}}
      e^{5 6}_{6 5},
      e^{5 7}_{7 5},
      e^{6 7}_{7 6}
      \\
      e^{6 5}_{5 6},
      e^{7 5}_{5 7},
      e^{7 6}_{6 7}
    \end{array}
  \right\},
  \\
  &&
  \hspace{-45pt}
  q
  \left\{
    \begin{array}{@{\hspace{0mm}}c@{\hspace{0mm}}}
      e^{2 5}_{5 2},
      e^{2 6}_{6 2},
      e^{3 5}_{5 3},
      e^{3 7}_{7 3},
      e^{4 6}_{6 4},
      e^{4 7}_{7 4}
      \\
      e^{5 2}_{2 5},
      e^{6 2}_{2 6},
      e^{5 3}_{3 5},
      e^{7 3}_{3 7},
      e^{6 4}_{4 6},
      e^{7 4}_{4 7}
    \end{array}
  \right\},
  \qquad
  q^{2}
  \left\{
    \begin{array}{@{\hspace{0mm}}c@{\hspace{0mm}}}
      e^{2 7}_{7 2},
      e^{4 5}_{5 4},
      e^{3 6}_{6 3}
      \\
      e^{7 2}_{2 7},
      e^{5 4}_{4 5},
      e^{6 3}_{3 6}
    \end{array}
  \right\},
  \qquad
  1
  \left\{
    \begin{array}{@{\hspace{0mm}}c@{\hspace{0mm}}}
      e^{1 8}_{8 1}
      \\
      e^{8 1}_{1 8}
    \end{array}
  \right\},
  \\[1mm]
  \hline
  \\[1mm]
  &&
  \hspace{-45pt}
  -
  \Delta
  q^{2}
  \left\{
    +
    1
    \left\{
    \begin{array}{@{\hspace{0mm}}c@{\hspace{0mm}}}
      e^{4 5}_{6 3}
      \\
      e^{6 3}_{4 5}
    \end{array}
    \right\},
    -
    q
    \left\{
    \begin{array}{@{\hspace{0mm}}c@{\hspace{0mm}}}
      e^{4 5}_{7 2}
      \\
      e^{7 2}_{4 5}
    \end{array}
    \right\},
    +
    1
    \left\{
    \begin{array}{@{\hspace{0mm}}c@{\hspace{0mm}}}
      e^{3 6}_{7 2}
      \\
      e^{7 2}_{3 6}
    \end{array}
    \right\}
  \right\},
  \\
  &&
  \hspace{-45pt}
  \Delta
  q^{3}
  [\alpha + 1]
  \left\{
    +
    \overline{q}
    \left\{
    \begin{array}{@{\hspace{0mm}}c@{\hspace{0mm}}}
      e^{5 4}_{6 3}
      \\
      e^{6 3}_{5 4}
    \end{array}
    \right\},
    -
    1
    \left\{
    \begin{array}{@{\hspace{0mm}}c@{\hspace{0mm}}}
      e^{5 4}_{7 2}
      \\
      e^{7 2}_{5 4}
    \end{array}
    \right\},
    +
    q
    \left\{
    \begin{array}{@{\hspace{0mm}}c@{\hspace{0mm}}}
      e^{6 3}_{7 2}
      \\
      e^{7 2}_{6 3}
    \end{array}
    \right\}
  \right\},
  \\
  &&
  \hspace{-45pt}
  \Delta
  q^{2}
  [\alpha]^{\frac{1}{2}}
  [\alpha + 2]^{\frac{1}{2}}
  \left\{
    -
    \overline{q}
    \left\{
    \begin{array}{@{\hspace{0mm}}c@{\hspace{0mm}}}
      e^{2 7}_{8 1}
      \\
      e^{8 1}_{2 7}
    \end{array}
    \right\},
    +
    1
    \left\{
    \begin{array}{@{\hspace{0mm}}c@{\hspace{0mm}}}
      e^{3 6}_{8 1}
      \\
      e^{8 1}_{3 6}
    \end{array}
    \right\},
    -
    q
    \left\{
    \begin{array}{@{\hspace{0mm}}c@{\hspace{0mm}}}
      e^{4 5}_{8 1}
      \\
      e^{8 1}_{4 5}
    \end{array}
    \right\}
  \right\},
  \\
  &&
  \hspace{-45pt}
  \Delta
  q^{- \alpha + 1}
  [\alpha]^{\frac{1}{2}}
  [\alpha + 1]^{\frac{1}{2}}
  \left\{
    -
    \overline{q}^{\frac{1}{2}}
    \left\{
    \begin{array}{@{\hspace{0mm}}c@{\hspace{0mm}}}
      e^{2 3}_{5 1},
      e^{2 4}_{6 1},
      e^{3 4}_{7 1}
      \\
      e^{5 1}_{2 3},
      e^{6 1}_{2 4},
      e^{7 1}_{3 4}
    \end{array}
    \right\},
    +
    q^{\frac{1}{2}}
    \left\{
    \begin{array}{@{\hspace{0mm}}c@{\hspace{0mm}}}
      e^{3 2}_{5 1},
      e^{4 2}_{6 1},
      e^{4 3}_{7 1}
      \\
      e^{5 1}_{3 2},
      e^{6 1}_{4 2},
      e^{7 1}_{4 3}
    \end{array}
    \right\}
  \right\},
  \\
  &&
  \hspace{-45pt}
  \Delta
  q^{\alpha + 3}
  [\alpha + 1]^{\frac{1}{2}}
  [\alpha + 2]^{\frac{1}{2}}
  \left\{
    -
    \overline{q}^{\frac{1}{2}}
    \left\{
    \begin{array}{@{\hspace{0mm}}c@{\hspace{0mm}}}
      e^{5 6}_{8 2},
      e^{5 7}_{8 3},
      e^{6 7}_{8 4}
      \\
      e^{8 2}_{5 6},
      e^{8 3}_{5 7},
      e^{8 4}_{6 7}
    \end{array}
    \right\},
    +
    q^{\frac{1}{2}}
    \left\{
    \begin{array}{@{\hspace{0mm}}c@{\hspace{0mm}}}
      e^{6 5}_{8 2},
      e^{7 5}_{8 3},
      e^{7 6}_{8 4}
      \\
      e^{8 2}_{6 5},
      e^{8 3}_{7 5},
      e^{8 4}_{7 6}
    \end{array}
    \right\}
  \right\},
  \\
  &&
  \hspace{-45pt}
  \Delta^2
  q^{3}
  [\alpha]^{\frac{1}{2}}
  [\alpha + 1]
  [\alpha + 2]^{\frac{1}{2}}
  \left\{
    +
    \overline{q}
    \left\{
    \begin{array}{@{\hspace{0mm}}c@{\hspace{0mm}}}
      e^{5 4}_{8 1}
      \\
      e^{8 1}_{5 4}
    \end{array}
    \right\},
    -
    1
    \left\{
    \begin{array}{@{\hspace{0mm}}c@{\hspace{0mm}}}
      e^{6 3}_{8 1}
      \\
      e^{8 1}_{6 3}
    \end{array}
    \right\},
    +
    q
    \left\{
    \begin{array}{@{\hspace{0mm}}c@{\hspace{0mm}}}
      e^{7 2}_{8 1}
      \\
      e^{8 1}_{7 2}
    \end{array}
    \right\}
  \right\},
\end{eqnarray*}
and the left handle $C$ has $8$ diagonal components:
\begin{eqnarray*}
  C
  =
  q^{- 3 \alpha - 3}
  \left\{
    \begin{array}{l}
      +
      q^{3}
      \left\{
                         e^{1}_{1}
      \right\}
      \\
      -
      q
      \left\{
        q^{2}            e^{2}_{2}, \; \:
                         e^{3}_{3}, \; \:
        \overline{q}^{2} e^{4}_{4}
      \right\}
      \\
      +
      \overline{q}
      \left\{
        \overline{q}^{2} e^{7}_{7}, \; \:
                         e^{6}_{6}, \; \:
        q^{2}            e^{5}_{5}
      \right\}
      \\
      -
      \overline{q}^{3}
      \left\{
                         e^{8}_{8}
      \right\}
    \end{array}
  \right\},
\end{eqnarray*}
using the scaling factors:
\begin{eqnarray*}
  \kappa_{\sigma}
  =
  q^{-3 \alpha},
  \qquad \qquad
  \kappa_{C}
  =
  q^{-3 \alpha}.
\end{eqnarray*}

\subsection*{Parameters for $LG^{4}$}

The reader will have by now appreciated the recurring patterns
in the components of our R matrices. To save space, we introduce
a little more notation, which eliminates the $q$ brackets altogether.
To whit, we write:
\begin{eqnarray*}
  A^z_i
  & \triangleq &
  {[ \alpha + i ]_q}^z,
  \qquad
  \mathrm{where~} z\in \{{\textstyle \frac{1}{2}}, 1\},
\end{eqnarray*}
and $i \in \{ 0, 1, 2, 3\}$.  With this notation, the braid
generator $\sigma$ has $758$ nonzero components:
\begin{eqnarray*}
  &&
  \hspace{-45pt}
  q^{- 4 \alpha}
  \left\{
    \begin{array}{@{\hspace{0mm}}c@{\hspace{0mm}}}
      e^{ 1, 1}_{ 1, 1}
    \end{array}
  \right\},
  \qquad
  q^{2}
  \left\{
    \begin{array}{@{\hspace{0mm}}c@{\hspace{0mm}}}
      e^{ 6, 6}_{ 6, 6},
      e^{ 7, 7}_{ 7, 7},
      e^{ 8, 8}_{ 8, 8},
      e^{ 9, 9}_{ 9, 9},
      e^{10,10}_{10,10},
      e^{11,11}_{11,11}
    \end{array}
  \right\},
  \qquad \hspace{-1pt}
  q^{4 \alpha + 12}
  \left\{
    \begin{array}{@{\hspace{0mm}}c@{\hspace{0mm}}}
      e^{16,16}_{16,16}
    \end{array}
  \right\},
  \\
  &&
  \hspace{-45pt}
  -
  q^{- 2 \alpha}
  \left\{
    \begin{array}{@{\hspace{0mm}}c@{\hspace{0mm}}}
      e^{ 2, 2}_{ 2, 2},
      e^{ 3, 3}_{ 3, 3},
      e^{ 4, 4}_{ 4, 4},
      e^{ 5, 5}_{ 5, 5}
    \end{array}
  \right\},
  \qquad
  -
  q^{2 \alpha + 6}
  \left\{
    \begin{array}{@{\hspace{0mm}}c@{\hspace{0mm}}}
      e^{15,15}_{15,15},
      e^{14,14}_{14,14},
      e^{13,13}_{13,13},
      e^{12,12}_{12,12}
    \end{array}
  \right\},
  \\[1mm]
  \hline
  \\[1mm]
  &&
  \hspace{-45pt}
  -
  \Delta
  q^{- 3 \alpha}
  A_{0}
  \left\{
    \begin{array}{@{\hspace{0mm}}c@{\hspace{0mm}}}
      e^{ 2, 1}_{ 2, 1},
      e^{ 3, 1}_{ 3, 1},
      e^{ 4, 1}_{ 4, 1},
      e^{ 5, 1}_{ 5, 1}
    \end{array}
  \right\},
  \qquad \hspace{-15pt}
  \Delta
  q^{3 \alpha + 9}
  A_{3}
  \left\{
    \begin{array}{@{\hspace{0mm}}c@{\hspace{0mm}}}
      e^{16,12}_{16,12},
      e^{16,13}_{16,13},
      e^{16,14}_{16,14},
      e^{16,15}_{16,15}
    \end{array}
  \right\},
  \\
  &&
  \hspace{-45pt}
  -
  \Delta
  q^{3}
  \left\{
    \begin{array}{@{\hspace{0mm}}c@{\hspace{0mm}}}
      e^{ 7, 6}_{ 7, 6},
      e^{ 8, 6}_{ 8, 6},
      e^{ 8, 7}_{ 8, 7},
      e^{ 9, 6}_{ 9, 6},
      e^{ 9, 7}_{ 9, 7},
      e^{10, 6}_{10, 6},
      e^{10, 8}_{10, 8},
      e^{10, 9}_{10, 9},
      e^{11, 7}_{11, 7},
      e^{11, 8}_{11, 8},
      e^{11, 9}_{11, 9},
      e^{11,10}_{11,10}
    \end{array}
  \right\},
  \\
  &&
  \hspace{-45pt}
  \Delta
  q^{- 2 \alpha + 1}
  \left\{
    \begin{array}{@{\hspace{0mm}}c@{\hspace{0mm}}}
      e^{ 3, 2}_{ 3, 2},
      e^{ 4, 2}_{ 4, 2},
      e^{ 4, 3}_{ 4, 3},
      e^{ 5, 2}_{ 5, 2},
      e^{ 5, 3}_{ 5, 3},
      e^{ 5, 4}_{ 5, 4}
    \end{array}
  \right\},
  \\
  &&
  \hspace{-45pt}
  \Delta
  q^{2 \alpha + 7}
  \left\{
    \begin{array}{@{\hspace{0mm}}c@{\hspace{0mm}}}
      e^{13,12}_{13,12},
      e^{14,12}_{14,12},
      e^{14,13}_{14,13},
      e^{15,12}_{15,12},
      e^{15,13}_{15,13},
      e^{15,14}_{15,14}
    \end{array}
  \right\},
  \\
  &&
  \hspace{-45pt}
  \Delta
  q^{- \alpha + 1}
  A_{1}
  \left\{
    \begin{array}{@{\hspace{0mm}}c@{\hspace{0mm}}}
      e^{ 6, 2}_{ 6, 2},
      e^{ 6, 3}_{ 6, 3},
      e^{ 7, 2}_{ 7, 2},
      e^{ 7, 4}_{ 7, 4},
      e^{ 8, 2}_{ 8, 2},
      e^{ 8, 5}_{ 8, 5},
      e^{ 9, 3}_{ 9, 3},
      e^{ 9, 4}_{ 9, 4},
      e^{10, 3}_{10, 3},
      e^{10, 5}_{10, 5},
      e^{11, 4}_{11, 4},
      e^{11, 5}_{11, 5}
    \end{array}
  \right\},
  \\
  &&
  \hspace{-45pt}
  \Delta
  q^{- \alpha + 3}
  (\overline{q}^2 - q^2)
  A_{1}
  \left\{
    \begin{array}{@{\hspace{0mm}}c@{\hspace{0mm}}}
      e^{ 9, 2}_{ 9, 2},
      e^{10, 2}_{10, 2},
      e^{11, 2}_{11, 2},
      e^{11, 3}_{11, 3}
    \end{array}
  \right\},
  \\
  &&
  \hspace{-45pt}
  \Delta^2
  q^{- 2 \alpha + 1}
  A_{0}
  A_{1}
  \left\{
    \begin{array}{@{\hspace{0mm}}c@{\hspace{0mm}}}
      e^{ 6, 1}_{ 6, 1},
      e^{ 7, 1}_{ 7, 1},
      e^{ 8, 1}_{ 8, 1},
      e^{ 9, 1}_{ 9, 1},
      e^{10, 1}_{10, 1},
      e^{11, 1}_{11, 1}
    \end{array}
  \right\},
  \\
  &&
  \hspace{-45pt}
  -
  \Delta
  q^{\alpha + 4}
  A_{2}
  \left\{
    \begin{array}{@{\hspace{0mm}}c@{\hspace{0mm}}}
      e^{12, 6}_{12, 6},
      e^{12, 7}_{12, 7},
      e^{12, 9}_{12, 9},
      e^{13, 6}_{13, 6},
      e^{13, 8}_{13, 8},
      e^{13,10}_{13,10}
      \\
      e^{14, 7}_{14, 7},
      e^{14, 8}_{14, 8},
      e^{14,11}_{14,11},
      e^{15, 9}_{15, 9},
      e^{15,10}_{15,10},
      e^{15,11}_{15,11}
    \end{array}
  \right\},
  \\
  &&
  \hspace{-45pt}
  -
  \Delta^2
  q^{- \alpha + 2}
  A_{1}
  \left\{
    \begin{array}{@{\hspace{0mm}}c@{\hspace{0mm}}}
      e^{ 7, 3}_{ 7, 3},
      e^{ 8, 3}_{ 8, 3},
      e^{ 8, 4}_{ 8, 4},
      e^{10, 4}_{10, 4}
    \end{array}
  \right\},
  \qquad
  \Delta^2
  q^{\alpha + 5}
  A_{2}
  \left\{
    \begin{array}{@{\hspace{0mm}}c@{\hspace{0mm}}}
      e^{13, 7}_{13, 7},
      e^{13, 9}_{13, 9},
      e^{14, 9}_{14, 9},
      e^{14,10}_{14,10}
    \end{array}
  \right\},
  \\
  &&
  \hspace{-45pt}
  -
  \Delta
  q^{\alpha + 6}
  (\overline{q}^2 - q^2)
  A_{2}
  \left\{
    \begin{array}{@{\hspace{0mm}}c@{\hspace{0mm}}}
      e^{14, 6}_{14, 6},
      e^{15, 6}_{15, 6},
      e^{15, 7}_{15, 7},
      e^{15, 8}_{15, 8}
    \end{array}
  \right\},
  \\
  &&
  \hspace{-45pt}
  -
  \Delta^2
  q^{3}
  A_{1}
  A_{2}
  \left\{
    \begin{array}{@{\hspace{0mm}}c@{\hspace{0mm}}}
      e^{12, 2}_{12, 2},
      e^{12, 3}_{12, 3},
      e^{12, 4}_{12, 4},
      e^{13, 2}_{13, 2},
      e^{13, 3}_{13, 3},
      e^{13, 5}_{13, 5}
      \\
      e^{14, 2}_{14, 2},
      e^{14, 4}_{14, 4},
      e^{14, 5}_{14, 5},
      e^{15, 3}_{15, 3},
      e^{15, 4}_{15, 4},
      e^{15, 5}_{15, 5}
    \end{array}
  \right\},
  \\
  &&
  \hspace{-45pt}
  -
  \Delta^3
  q^{- \alpha + 3}
  A_{0}
  A_{1}
  A_{2}
  \left\{
    \begin{array}{@{\hspace{0mm}}c@{\hspace{0mm}}}
      e^{12, 1}_{12, 1},
      e^{13, 1}_{13, 1},
      e^{14, 1}_{14, 1},
      e^{15, 1}_{15, 1}
    \end{array}
  \right\},
  \\
  &&
  \hspace{-45pt}
  \Delta^2
  q^{2 \alpha + 7}
  A_{2}
  A_{3}
  \left\{
    \begin{array}{@{\hspace{0mm}}c@{\hspace{0mm}}}
      e^{16, 6}_{16, 6},
      e^{16, 7}_{16, 7},
      e^{16, 8}_{16, 8},
      e^{16, 9}_{16, 9},
      e^{16,10}_{16,10},
      e^{16,11}_{16,11}
    \end{array}
  \right\},
  \\
  &&
  \hspace{-45pt}
  \Delta^3
  q^{\alpha + 6}
  A_{1}
  A_{2}
  A_{3}
  \left\{
    \begin{array}{@{\hspace{0mm}}c@{\hspace{0mm}}}
      e^{16, 2}_{16, 2},
      e^{16, 3}_{16, 3},
      e^{16, 4}_{16, 4},
      e^{16, 5}_{16, 5}
    \end{array}
  \right\},
  \qquad
  \Delta^3
  q^{6}
  (\overline{q}^{2} + 1 + q^{2})
  A_{1}
  A_{2}
  \left\{
    \begin{array}{@{\hspace{0mm}}c@{\hspace{0mm}}}
      e^{15, 2}_{15, 2}
    \end{array}
  \right\},
  \\
  &&
  \hspace{-45pt}
  \Delta^2
  q^{4}
  \left\{
    \begin{array}{@{\hspace{0mm}}c@{\hspace{0mm}}}
      e^{10, 7}_{10, 7}
    \end{array}
  \right\},
  \qquad
  \Delta^2
  q^{5}
  ( q + \overline{q} )
  \left\{
    \begin{array}{@{\hspace{0mm}}c@{\hspace{0mm}}}
      e^{11, 6}_{11, 6}
    \end{array}
  \right\},
  \qquad
  \Delta^3
  q^{4}
  A_{1}
  A_{2}
  \left\{
    \begin{array}{@{\hspace{0mm}}c@{\hspace{0mm}}}
      e^{13, 4}_{13, 4}
    \end{array}
  \right\},
  \\
  &&
  \hspace{-45pt}
  -
  \Delta^2
  q^{5}
  (\overline{q}^2 - q^2)
  A_{1}
  A_{2}
  \left\{
    \begin{array}{@{\hspace{0mm}}c@{\hspace{0mm}}}
      e^{14, 3}_{14, 3}
    \end{array}
  \right\},
  \qquad \hspace{-9pt}
  \Delta^4
  q^{6}
  A_{0}
  A_{1}
  A_{2}
  A_{3}
  \left\{
    \begin{array}{@{\hspace{0mm}}c@{\hspace{0mm}}}
      e^{16, 1}_{16, 1}
    \end{array}
  \right\},
  \\[1mm]
  \hline
  \\[1mm]
  &&
  \hspace{-45pt}
  q^{- 3 \alpha}
  \left\{
    \begin{array}{@{\hspace{0mm}}c@{\hspace{0mm}}}
      e^{ 1, 2}_{ 2, 1},
      e^{ 1, 3}_{ 3, 1},
      e^{ 1, 4}_{ 4, 1},
      e^{ 1, 5}_{ 5, 1}
      \\
      e^{ 2, 1}_{ 1, 2},
      e^{ 3, 1}_{ 1, 3},
      e^{ 4, 1}_{ 1, 4},
      e^{ 5, 1}_{ 1, 5}
    \end{array}
  \right\},
  \qquad
  q^{2 \alpha + 6}
  \left\{
    \begin{array}{@{\hspace{0mm}}c@{\hspace{0mm}}}
      e^{11,16}_{16,11},
      e^{10,16}_{16,10},
      e^{ 9,16}_{16, 9},
      e^{ 8,16}_{16, 8},
      e^{ 7,16}_{16, 7},
      e^{ 6,16}_{16, 6}
      \\
      e^{16,11}_{11,16},
      e^{16,10}_{10,16},
      e^{16, 9}_{ 9,16},
      e^{16, 8}_{ 8,16},
      e^{16, 7}_{ 7,16},
      e^{16, 6}_{ 6,16}
    \end{array}
  \right\},
  \\
  &&
  \hspace{-45pt}
  q^{3 \alpha + 9}
  \left\{
    \begin{array}{@{\hspace{0mm}}c@{\hspace{0mm}}}
      e^{16,15}_{15,16},
      e^{16,14}_{14,16},
      e^{16,13}_{13,16},
      e^{16,12}_{12,16}
      \\
      e^{15,16}_{16,15},
      e^{14,16}_{16,14},
      e^{13,16}_{16,13},
      e^{12,16}_{16,12}
    \end{array}
  \right\},
  \qquad \hspace{-3pt}
  q^{- 2 \alpha}
  \left\{
    \begin{array}{@{\hspace{0mm}}c@{\hspace{0mm}}}
      e^{ 1, 6}_{ 6, 1},
      e^{ 1, 7}_{ 7, 1},
      e^{ 1, 8}_{ 8, 1},
      e^{ 1, 9}_{ 9, 1},
      e^{ 1,10}_{10, 1},
      e^{ 1,11}_{11, 1}
      \\
      e^{ 6, 1}_{ 1, 6},
      e^{ 7, 1}_{ 1, 7},
      e^{ 8, 1}_{ 1, 8},
      e^{ 9, 1}_{ 1, 9},
      e^{10, 1}_{ 1,10},
      e^{11, 1}_{ 1,11}
    \end{array}
  \right\},
  \\
  &&
  \hspace{-45pt}
  q^{- \alpha}
  \left\{
    \begin{array}{@{\hspace{0mm}}c@{\hspace{0mm}}}
      e^{ 1,12}_{12, 1},
      e^{ 1,13}_{13, 1},
      e^{ 1,14}_{14, 1},
      e^{ 1,15}_{15, 1}
      \\
      e^{12, 1}_{ 1,12},
      e^{13, 1}_{ 1,13},
      e^{14, 1}_{ 1,14},
      e^{15, 1}_{ 1,15}
    \end{array}
  \right\},
  \qquad
  -
  q^{- 2 \alpha + 1}
  \left\{
    \begin{array}{@{\hspace{0mm}}c@{\hspace{0mm}}}
      e^{ 2, 3}_{ 3, 2},
      e^{ 2, 4}_{ 4, 2},
      e^{ 2, 5}_{ 5, 2},
      e^{ 3, 4}_{ 4, 3},
      e^{ 3, 5}_{ 5, 3},
      e^{ 4, 5}_{ 5, 4}
      \\
      e^{ 3, 2}_{ 2, 3},
      e^{ 4, 2}_{ 2, 4},
      e^{ 5, 2}_{ 2, 5},
      e^{ 4, 3}_{ 3, 4},
      e^{ 5, 3}_{ 3, 5},
      e^{ 5, 4}_{ 4, 5}
    \end{array}
  \right\},
  \\
  &&
  \hspace{-45pt}
  q^{\alpha + 3}
  \left\{
    \begin{array}{@{\hspace{0mm}}c@{\hspace{0mm}}}
      e^{ 5,16}_{16, 5},
      e^{ 4,16}_{16, 4},
      e^{ 3,16}_{16, 3},
      e^{ 2,16}_{16, 2}
      \\
      e^{16, 5}_{ 5,16},
      e^{16, 4}_{ 4,16},
      e^{16, 3}_{ 3,16},
      e^{16, 2}_{ 2,16}
    \end{array}
  \right\},
  \quad
  \hspace{-4pt}
  -
  q^{2 \alpha + 7}
  \left\{
    \begin{array}{@{\hspace{0mm}}c@{\hspace{0mm}}}
      e^{15,14}_{14,15},
      e^{15,13}_{13,15},
      e^{15,12}_{12,15},
      e^{14,13}_{13,14},
      e^{14,12}_{12,14},
      e^{13,12}_{12,13}
      \\
      e^{14,15}_{15,14},
      e^{13,15}_{15,13},
      e^{12,15}_{15,12},
      e^{13,14}_{14,13},
      e^{12,14}_{14,12},
      e^{12,13}_{13,12}
    \end{array}
  \right\},
  \\
  &&
  \hspace{-45pt}
  q^{- \alpha + 1}
  \left\{
    \begin{array}{@{\hspace{0mm}}c@{\hspace{0mm}}}
      e^{ 2, 6}_{ 6, 2},
      e^{ 2, 7}_{ 7, 2},
      e^{ 2, 8}_{ 8, 2},
      e^{ 3, 6}_{ 6, 3},
      e^{ 3, 9}_{ 9, 3},
      e^{ 3,10}_{10, 3},
      e^{ 4, 7}_{ 7, 4},
      e^{ 4, 9}_{ 9, 4},
      e^{ 4,11}_{11, 4},
      e^{ 5, 8}_{ 8, 5},
      e^{ 5,10}_{10, 5},
      e^{ 5,11}_{11, 5}
      \\
      e^{ 6, 2}_{ 2, 6},
      e^{ 7, 2}_{ 2, 7},
      e^{ 8, 2}_{ 2, 8},
      e^{ 6, 3}_{ 3, 6},
      e^{ 9, 3}_{ 3, 9},
      e^{10, 3}_{ 3,10},
      e^{ 7, 4}_{ 4, 7},
      e^{ 9, 4}_{ 4, 9},
      e^{11, 4}_{ 4,11},
      e^{ 8, 5}_{ 5, 8},
      e^{10, 5}_{ 5,10},
      e^{11, 5}_{ 5,11}
    \end{array}
  \right\},
  \\
  &&
  \hspace{-45pt}
  q^{\alpha + 4}
  \left\{
    \begin{array}{@{\hspace{0mm}}c@{\hspace{0mm}}}
      e^{15,10}_{10,15},
      e^{15,11}_{11,15},
      e^{15, 9}_{ 9,15},
      e^{14,11}_{11,14},
      e^{14, 8}_{ 8,14},
      e^{14, 7}_{ 7,14},
      e^{13,10}_{10,13},
      e^{13, 8}_{ 8,13},
      e^{13, 6}_{ 6,13},
      e^{12, 9}_{ 9,12},
      e^{12, 7}_{ 7,12},
      e^{12, 6}_{ 6,12}
      \\
      e^{11,15}_{15,11},
      e^{10,15}_{15,10},
      e^{ 9,15}_{15, 9},
      e^{11,14}_{14,11},
      e^{ 8,14}_{14, 8},
      e^{ 7,14}_{14, 7},
      e^{10,13}_{13,10},
      e^{ 8,13}_{13, 8},
      e^{ 6,13}_{13, 6},
      e^{ 9,12}_{12, 9},
      e^{ 7,12}_{12, 7},
      e^{ 6,12}_{12, 6}
    \end{array}
  \right\},
  \\
  &&
  \hspace{-45pt}
  q^{- \alpha + 2}
  \left\{
    \begin{array}{@{\hspace{0mm}}c@{\hspace{0mm}}}
      e^{ 2, 9}_{ 9, 2},
      e^{ 2,10}_{10, 2},
      e^{ 2,11}_{11, 2},
      e^{ 3, 7}_{ 7, 3},
      e^{ 3, 8}_{ 8, 3},
      e^{ 3,11}_{11, 3},
      e^{ 4, 6}_{ 6, 4},
      e^{ 4, 8}_{ 8, 4},
      e^{ 4,10}_{10, 4},
      e^{ 5, 6}_{ 6, 5},
      e^{ 5, 7}_{ 7, 5},
      e^{ 5, 9}_{ 9, 5}
      \\
      e^{ 9, 2}_{ 2, 9},
      e^{10, 2}_{ 2,10},
      e^{11, 2}_{ 2,11},
      e^{ 7, 3}_{ 3, 7},
      e^{ 8, 3}_{ 3, 8},
      e^{11, 3}_{ 3,11},
      e^{ 6, 4}_{ 4, 6},
      e^{ 8, 4}_{ 4, 8},
      e^{10, 4}_{ 4,10},
      e^{ 6, 5}_{ 5, 6},
      e^{ 7, 5}_{ 5, 7},
      e^{ 9, 5}_{ 5, 9}
    \end{array}
  \right\},
  \\
  &&
  \hspace{-45pt}
  q^{3}
  \left\{
    \begin{array}{@{\hspace{0mm}}c@{\hspace{0mm}}}
        e^{ 6, 7}_{ 7, 6},
        e^{ 6, 8}_{ 8, 6},
        e^{ 6, 9}_{ 9, 6},
        e^{ 6,10}_{10, 6},
        e^{ 7, 8}_{ 8, 7},
        e^{ 7, 9}_{ 9, 7},
        e^{ 7,11}_{11, 7},
        e^{ 8,10}_{10, 8},
        e^{ 8,11}_{11, 8},
        e^{ 9,10}_{10, 9},
        e^{ 9,11}_{11, 9},
        e^{10,11}_{11,10}
      \\
        e^{ 7, 6}_{ 6, 7},
        e^{ 8, 6}_{ 6, 8},
        e^{ 9, 6}_{ 6, 9},
        e^{10, 6}_{ 6,10},
        e^{ 8, 7}_{ 7, 8},
        e^{ 9, 7}_{ 7, 9},
        e^{10, 8}_{ 8,10},
        e^{11, 7}_{ 7,11},
        e^{11, 8}_{ 8,11},
        e^{10, 9}_{ 9,10},
        e^{11, 9}_{ 9,11},
        e^{11,10}_{10,11}
    \end{array}
  \right\},
  \\
  &&
  \hspace{-45pt}
  1
  \left\{
    \begin{array}{@{\hspace{0mm}}c@{\hspace{0mm}}}
      e^{ 1,16}_{16, 1}
      \\
      e^{16, 1}_{ 1,16}
    \end{array}
  \right\},
  \qquad
  -
  q^{3}
  \left\{
    \begin{array}{@{\hspace{0mm}}c@{\hspace{0mm}}}
      e^{ 2,15}_{15, 2},
      e^{ 3,14}_{14, 3},
      e^{ 4,13}_{13, 4},
      e^{ 5,12}_{12, 5}
      \\
      e^{12, 5}_{ 5,12},
      e^{13, 4}_{ 4,13},
      e^{14, 3}_{ 3,14},
      e^{15, 2}_{ 2,15}
    \end{array}
  \right\},
  \qquad
  q^{4}
  \left\{
    \begin{array}{@{\hspace{0mm}}c@{\hspace{0mm}}}
      e^{ 6,11}_{11, 6},
      e^{ 7,10}_{10, 7},
      e^{ 8, 9}_{ 9, 8}
      \\
      e^{11, 6}_{ 6,11},
      e^{10, 7}_{ 7,10},
      e^{ 9, 8}_{ 8, 9}
    \end{array}
  \right\},
  \\
  &&
  \hspace{-45pt}
  q^{\alpha + 5}
  \left\{
    \begin{array}{@{\hspace{0mm}}c@{\hspace{0mm}}}
      e^{ 6,14}_{14, 6},
      e^{ 6,15}_{15, 6},
      e^{ 7,13}_{13, 7},
      e^{ 7,15}_{15, 7},
      e^{ 8,12}_{12, 8},
      e^{ 8,15}_{15, 8},
      e^{ 9,13}_{13, 9},
      e^{ 9,14}_{14, 9},
      e^{10,12}_{12,10},
      e^{10,14}_{14,10},
      e^{11,12}_{12,11},
      e^{11,13}_{13,11}
      \\
      e^{14, 6}_{ 6,14},
      e^{15, 6}_{ 6,15},
      e^{13, 7}_{ 7,13},
      e^{15, 7}_{ 7,15},
      e^{12, 8}_{ 8,12},
      e^{15, 8}_{ 8,15},
      e^{13, 9}_{ 9,13},
      e^{14, 9}_{ 9,14},
      e^{12,10}_{10,12},
      e^{14,10}_{10,14},
      e^{12,11}_{11,12},
      e^{13,11}_{11,13}
    \end{array}
  \right\},
  \\
  &&
  \hspace{-45pt}
  -
  q^{2}
  \left\{
    \begin{array}{@{\hspace{0mm}}c@{\hspace{0mm}}}
      e^{ 2,12}_{12, 2},
      e^{ 2,13}_{13, 2},
      e^{ 2,14}_{14, 2},
      e^{ 3,12}_{12, 3},
      e^{ 3,13}_{13, 3},
      e^{ 3,15}_{15, 3},
      e^{ 4,12}_{12, 4},
      e^{ 4,14}_{14, 4},
      e^{ 4,15}_{15, 4},
      e^{ 5,13}_{13, 5},
      e^{ 5,14}_{14, 5},
      e^{ 5,15}_{15, 5}
      \\
      e^{12, 2}_{ 2,12},
      e^{12, 3}_{ 3,12},
      e^{12, 4}_{ 4,12},
      e^{13, 2}_{ 2,13},
      e^{13, 3}_{ 3,13},
      e^{13, 5}_{ 5,13},
      e^{14, 2}_{ 2,14},
      e^{14, 4}_{ 4,14},
      e^{14, 5}_{ 5,14},
      e^{15, 3}_{ 3,15},
      e^{15, 4}_{ 4,15},
      e^{15, 5}_{ 5,15}
    \end{array}
  \right\},
  \\[1mm]
  \hline
  \\[1mm]
  &&
  \hspace{-45pt}
  \Delta
  q^{4}
  \left\{
    +
    \overline{q}
    \left\{
    \begin{array}{@{\hspace{0mm}}c@{\hspace{0mm}}}
      e^{ 3,14}_{15, 2},
      e^{ 4,13}_{14, 3},
      e^{ 5,12}_{13, 4}
      \\
      e^{15, 2}_{ 3,14},
      e^{14, 3}_{ 4,13},
      e^{13, 4}_{ 5,12}
    \end{array}
    \right\},
    -
    1
    \left\{
    \begin{array}{@{\hspace{0mm}}c@{\hspace{0mm}}}
      e^{ 4,13}_{15, 2},
      e^{ 5,12}_{14, 3}
      \\
      e^{15, 2}_{ 4,13},
      e^{14, 3}_{ 5,12}
    \end{array}
    \right\},
    +
    q
    \left\{
    \begin{array}{@{\hspace{0mm}}c@{\hspace{0mm}}}
      e^{ 5,12}_{15, 2}
      \\
      e^{15, 2}_{ 5,12}
    \end{array}
    \right\}
  \right\},
  \\
  &&
  \hspace{-45pt}
  \Delta
  q^{5}
  \left\{
    -
    \overline{q}
    \left\{
    \begin{array}{@{\hspace{0mm}}c@{\hspace{0mm}}}
      e^{ 7,10}_{11, 6},
      e^{ 8, 9}_{10, 7},
      e^{ 9, 8}_{10, 7}
      \\
      e^{11, 6}_{ 7,10},
      e^{10, 7}_{ 8, 9},
      e^{10, 7}_{ 9, 8}
    \end{array}
    \right\},
    +
    1
    \left\{
    \begin{array}{@{\hspace{0mm}}c@{\hspace{0mm}}}
      e^{ 8, 9}_{11, 6},
      e^{ 9, 8}_{11, 6}
      \\
      e^{11, 6}_{ 8, 9},
      e^{11, 6}_{ 9, 8}
    \end{array}
    \right\},
    -
    q
    \left\{
    \begin{array}{@{\hspace{0mm}}c@{\hspace{0mm}}}
      e^{10, 7}_{11, 6}
      \\
      e^{11, 6}_{10, 7}
    \end{array}
    \right\}
  \right\},
  \\
  &&
  \hspace{-45pt}
  -
  \Delta
  q^{- \alpha + 2}
  \left\{
    \begin{array}{@{\hspace{0mm}}c@{\hspace{0mm}}}
      e^{ 3, 7}_{ 9, 2},
      e^{ 3, 8}_{10, 2},
      e^{ 4, 6}_{ 7, 3},
      e^{ 4, 8}_{11, 2},
      e^{ 4,10}_{11, 3},
      e^{ 5, 6}_{ 8, 3},
      e^{ 5, 7}_{ 8, 4},
      e^{ 5, 9}_{10, 4}
      \\
      e^{ 9, 2}_{ 3, 7},
      e^{10, 2}_{ 3, 8},
      e^{ 7, 3}_{ 4, 6},
      e^{11, 2}_{ 4, 8},
      e^{11, 3}_{ 4,10},
      e^{ 8, 3}_{ 5, 6},
      e^{ 8, 4}_{ 5, 7},
      e^{10, 4}_{ 5, 9}
    \end{array}
  \right\},
  \\
  &&
  \hspace{-45pt}
  -
  \Delta
  q^{\alpha + 5}
  \left\{
    \begin{array}{@{\hspace{0mm}}c@{\hspace{0mm}}}
      e^{ 7,13}_{14, 6},
      e^{ 8,12}_{13, 7},
      e^{ 9,13}_{15, 6},
      e^{ 9,14}_{15, 7},
      e^{10,12}_{13, 9},
      e^{10,14}_{15, 8},
      e^{11,12}_{14, 9},
      e^{11,13}_{14,10}
      \\
      e^{14, 6}_{ 7,13},
      e^{13, 7}_{ 8,12},
      e^{15, 6}_{ 9,13},
      e^{15, 7}_{ 9,14},
      e^{13, 9}_{10,12},
      e^{15, 8}_{10,14},
      e^{14, 9}_{11,12},
      e^{14,10}_{11,13}
    \end{array}
  \right\},
  \\
  &&
  \hspace{-45pt}
  \Delta
  q^{- \alpha + 3}
  \left\{
    \begin{array}{@{\hspace{0mm}}c@{\hspace{0mm}}}
      e^{ 4, 6}_{ 9, 2},
      e^{ 5, 6}_{10, 2},
      e^{ 5, 7}_{11, 2},
      e^{ 5, 9}_{11, 3}
      \\
      e^{ 9, 2}_{ 4, 6},
      e^{10, 2}_{ 5, 6},
      e^{11, 2}_{ 5, 7},
      e^{11, 3}_{ 5, 9}
    \end{array}
  \right\},
  \qquad
  \Delta
  q^{\alpha + 6}
  \left\{
    \begin{array}{@{\hspace{0mm}}c@{\hspace{0mm}}}
      e^{ 8,12}_{14, 6},
      e^{10,12}_{15, 6},
      e^{11,12}_{15, 7},
      e^{11,13}_{15, 8}
      \\
      e^{14, 6}_{ 8,12},
      e^{15, 6}_{10,12},
      e^{15, 7}_{11,12},
      e^{15, 8}_{11,13}
    \end{array}
  \right\},
  \\
  &&
  \hspace{-45pt}
  -
  \Delta
  q^{- 2 \alpha + \frac{1}{2}}
  A_{0}^{\frac{1}{2}}
  A_{1}^{\frac{1}{2}}
  \left\{
    \begin{array}{@{\hspace{0mm}}c@{\hspace{0mm}}}
      e^{ 2, 3}_{ 6, 1},
      e^{ 2, 4}_{ 7, 1},
      e^{ 2, 5}_{ 8, 1},
      e^{ 3, 4}_{ 9, 1},
      e^{ 3, 5}_{10, 1},
      e^{ 4, 5}_{11, 1}
      \\
      e^{ 6, 1}_{ 2, 3},
      e^{ 7, 1}_{ 2, 4},
      e^{ 8, 1}_{ 2, 5},
      e^{ 9, 1}_{ 3, 4},
      e^{10, 1}_{ 3, 5},
      e^{11, 1}_{ 4, 5}
    \end{array}
  \right\},
  \\
  &&
  \hspace{-45pt}
  \Delta
  q^{- 2 \alpha + \frac{3}{2}}
  A_{0}^{\frac{1}{2}}
  A_{1}^{\frac{1}{2}}
  \left\{
    \begin{array}{@{\hspace{0mm}}c@{\hspace{0mm}}}
      e^{ 3, 2}_{ 6, 1},
      e^{ 4, 2}_{ 7, 1},
      e^{ 5, 2}_{ 8, 1},
      e^{ 4, 3}_{ 9, 1},
      e^{ 5, 3}_{10, 1},
      e^{ 5, 4}_{11, 1}
      \\
      e^{ 6, 1}_{ 3, 2},
      e^{ 7, 1}_{ 4, 2},
      e^{ 8, 1}_{ 5, 2},
      e^{ 9, 1}_{ 4, 3},
      e^{10, 1}_{ 5, 3},
      e^{11, 1}_{ 5, 4}
    \end{array}
  \right\},
  \\
  &&
  \hspace{-45pt}
  \Delta
  q^{- \alpha + 3}
  A_{1}
  \left\{
    \begin{array}{@{\hspace{0mm}}c@{\hspace{0mm}}}
      +
      \overline{q}
      \left\{
      \begin{array}{@{\hspace{0mm}}c@{\hspace{0mm}}}
        e^{ 6, 4}_{ 7, 3},
        e^{ 6, 5}_{ 8, 3},
        e^{ 7, 5}_{ 8, 4},
        e^{ 9, 5}_{10, 4}
        \\
        e^{ 7, 3}_{ 6, 4},
        e^{ 8, 3}_{ 6, 5},
        e^{ 8, 4}_{ 7, 5},
        e^{10, 4}_{ 9, 5}
      \end{array}
      \right\}
      \\
      -
      1
      \left\{
      \begin{array}{@{\hspace{0mm}}c@{\hspace{0mm}}}
        e^{ 6, 4}_{ 9, 2},
        e^{ 6, 5}_{10, 2},
        e^{ 7, 5}_{11, 2},
        e^{ 9, 5}_{11, 3}
        \\
        e^{ 9, 2}_{ 6, 4},
        e^{10, 2}_{ 6, 5},
        e^{11, 2}_{ 7, 5},
        e^{11, 3}_{ 9, 5}
      \end{array}
      \right\}
      \\
      +
      q
      \left\{
      \begin{array}{@{\hspace{0mm}}c@{\hspace{0mm}}}
        e^{ 7, 3}_{ 9, 2},
        e^{ 8, 3}_{10, 2},
        e^{ 8, 4}_{11, 2},
        e^{10, 4}_{11, 3}
        \\
        e^{ 9, 2}_{ 7, 3},
        e^{10, 2}_{ 8, 3},
        e^{11, 2}_{ 8, 4},
        e^{11, 3}_{10, 4}
      \end{array}
      \right\}
    \end{array}
  \right\},
  \\
  &&
  \hspace{-45pt}
  \Delta
  q^{\alpha + 6}
  A_{2}
  \left\{
    \begin{array}{@{\hspace{0mm}}c@{\hspace{0mm}}}
      -
      \overline{q}
      \left\{
      \begin{array}{@{\hspace{0mm}}c@{\hspace{0mm}}}
        e^{12, 8}_{13, 7},
        e^{12,10}_{13, 9},
        e^{12,11}_{14, 9},
        e^{13,11}_{14,10}
        \\
        e^{13, 7}_{12, 8},
        e^{13, 9}_{12,10},
        e^{14, 9}_{12,11},
        e^{14,10}_{13,11}
      \end{array}
      \right\}
      \\
      + 1
      \left\{
      \begin{array}{@{\hspace{0mm}}c@{\hspace{0mm}}}
        e^{12, 8}_{14, 6},
        e^{12,10}_{15, 6},
        e^{12,11}_{15, 7},
        e^{13,11}_{15, 8}
        \\
        e^{14, 6}_{12, 8},
        e^{15, 6}_{12,10},
        e^{15, 7}_{12,11},
        e^{15, 8}_{13,11}
      \end{array}
      \right\}
      \\
      -
      q
      \left\{
      \begin{array}{@{\hspace{0mm}}c@{\hspace{0mm}}}
        e^{13, 7}_{14, 6},
        e^{13, 9}_{15, 6},
        e^{14, 9}_{15, 7},
        e^{14,10}_{15, 8}
        \\
        e^{14, 6}_{13, 7},
        e^{15, 6}_{13, 9},
        e^{15, 7}_{14, 9},
        e^{15, 8}_{14,10}
      \end{array}
      \right\}
    \end{array}
  \right\},
  \\
  &&
  \hspace{-45pt}
  \Delta
  q^{- \alpha + 2}
  A_{0}^{\frac{1}{2}}
  A_{2}^{\frac{1}{2}}
  \left\{
    \begin{array}{@{\hspace{0mm}}c@{\hspace{0mm}}}
      -
      \overline{q}
      \left\{
      \begin{array}{@{\hspace{0mm}}c@{\hspace{0mm}}}
        e^{ 2, 9}_{12, 1},
        e^{ 2,10}_{13, 1},
        e^{ 2,11}_{14, 1},
        e^{ 3,11}_{15, 1}
        \\
        e^{12, 1}_{ 2, 9},
        e^{13, 1}_{ 2,10},
        e^{14, 1}_{ 2,11},
        e^{15, 1}_{ 3,11}
      \end{array}
      \right\}
      \\
      +
      1
      \left\{
      \begin{array}{@{\hspace{0mm}}c@{\hspace{0mm}}}
        e^{ 3, 7}_{12, 1},
        e^{ 3, 8}_{13, 1},
        e^{ 4, 8}_{14, 1},
        e^{ 4,10}_{15, 1}
        \\
        e^{12, 1}_{ 3, 7},
        e^{13, 1}_{ 3, 8},
        e^{14, 1}_{ 4, 8},
        e^{15, 1}_{ 4,10}
      \end{array}
      \right\}
      \\
      -
      q
      \left\{
      \begin{array}{@{\hspace{0mm}}c@{\hspace{0mm}}}
        e^{ 4, 6}_{12, 1},
        e^{ 5, 6}_{13, 1},
        e^{ 5, 7}_{14, 1},
        e^{ 5, 9}_{15, 1}
        \\
        e^{12, 1}_{ 4, 6},
        e^{13, 1}_{ 5, 6},
        e^{14, 1}_{ 5, 7},
        e^{15, 1}_{ 5, 9}
      \end{array}
      \right\}
    \end{array}
  \right\},
  \\
  & &
  \hspace{-45pt}
  \Delta
  q^{\alpha + 5}
  A_{1}^{\frac{1}{2}}
  A_{3}^{\frac{1}{2}}
  \left\{
    \begin{array}{@{\hspace{0mm}}c@{\hspace{0mm}}}
      +
      \overline{q}
      \left\{
      \begin{array}{@{\hspace{0mm}}c@{\hspace{0mm}}}
        e^{ 6,14}_{16, 2},
        e^{ 6,15}_{16, 3},
        e^{ 7,15}_{16, 4},
        e^{ 8,15}_{16, 5}
        \\
        e^{16, 2}_{ 6,14},
        e^{16, 3}_{ 6,15},
        e^{16, 4}_{ 7,15},
        e^{16, 5}_{ 8,15}
      \end{array}
      \right\}
      \\
      -
      1
      \left\{
      \begin{array}{@{\hspace{0mm}}c@{\hspace{0mm}}}
        e^{ 7,13}_{16, 2},
        e^{ 9,13}_{16, 3},
        e^{ 9,14}_{16, 4},
        e^{10,14}_{16, 5}
        \\
        e^{16, 2}_{ 7,13},
        e^{16, 3}_{ 9,13},
        e^{16, 4}_{ 9,14},
        e^{16, 5}_{10,14}
      \end{array}
      \right\}
      \\
      +
      q
      \left\{
      \begin{array}{@{\hspace{0mm}}c@{\hspace{0mm}}}
        e^{ 8,12}_{16, 2},
        e^{10,12}_{16, 3},
        e^{11,12}_{16, 4},
        e^{11,13}_{16, 5}
        \\
        e^{16, 2}_{ 8,12},
        e^{16, 3}_{10,12},
        e^{16, 4}_{11,12},
        e^{16, 5}_{11,13}
      \end{array}
      \right\}
    \end{array}
  \right\},
  \\
  &&
  \hspace{-45pt}
  \Delta^2
  q^{- \alpha + 3}
  A_{0}^{\frac{1}{2}}
  A_{1}
  A_{2}^{\frac{1}{2}}
  \left\{
    \begin{array}{@{\hspace{0mm}}c@{\hspace{0mm}}}
      +
      \overline{q}
      \left\{
      \begin{array}{@{\hspace{0mm}}c@{\hspace{0mm}}}
        e^{ 6, 4}_{12, 1},
        e^{ 6, 5}_{13, 1},
        e^{ 7, 5}_{14, 1},
        e^{ 9, 5}_{15, 1}
        \\
        e^{12, 1}_{ 6, 4},
        e^{13, 1}_{ 6, 5},
        e^{14, 1}_{ 7, 5},
        e^{15, 1}_{ 9, 5}
      \end{array}
      \right\}
      \\
      -
      1
      \left\{
      \begin{array}{@{\hspace{0mm}}c@{\hspace{0mm}}}
        e^{ 7, 3}_{12, 1},
        e^{ 8, 3}_{13, 1},
        e^{ 8, 4}_{14, 1},
        e^{10, 4}_{15, 1}
      \\
        e^{12, 1}_{ 7, 3},
        e^{13, 1}_{ 8, 3},
        e^{14, 1}_{ 8, 4},
        e^{15, 1}_{10, 4}
      \end{array}
      \right\}
      \\
      +
      q
      \left\{
      \begin{array}{@{\hspace{0mm}}c@{\hspace{0mm}}}
        e^{ 9, 2}_{12, 1},
        e^{10, 2}_{13, 1},
        e^{11, 2}_{14, 1},
        e^{11, 3}_{15, 1}
        \\
        e^{12, 1}_{ 9, 2},
        e^{13, 1}_{10, 2},
        e^{14, 1}_{11, 2},
        e^{15, 1}_{11, 3}
      \end{array}
      \right\}
    \end{array}
  \right\},
  \\
  &&
  \hspace{-45pt}
  \Delta^2
  q^{\alpha + 6}
  A_{1}^{\frac{1}{2}}
  A_{2}
  A_{3}^{\frac{1}{2}}
  \left\{
    \begin{array}{@{\hspace{0mm}}c@{\hspace{0mm}}}
      +
      \overline{q}
      \left\{
      \begin{array}{@{\hspace{0mm}}c@{\hspace{0mm}}}
        e^{12, 8}_{16, 2},
        e^{12,10}_{16, 3},
        e^{12,11}_{16, 4},
        e^{13,11}_{16, 5}
        \\
        e^{16, 2}_{12, 8},
        e^{16, 3}_{12,10},
        e^{16, 4}_{12,11},
        e^{16, 5}_{13,11}
      \end{array}
      \right\}
      \\
      - 1
      \left\{
      \begin{array}{@{\hspace{0mm}}c@{\hspace{0mm}}}
        e^{13, 7}_{16, 2},
        e^{13, 9}_{16, 3},
        e^{14, 9}_{16, 4},
        e^{14,10}_{16, 5}
        \\
        e^{16, 2}_{13, 7},
        e^{16, 3}_{13, 9},
        e^{16, 4}_{14, 9},
        e^{16, 5}_{14,10}
      \end{array}
      \right\}
      \\
      + q
      \left\{
      \begin{array}{@{\hspace{0mm}}c@{\hspace{0mm}}}
        e^{14, 6}_{16, 2},
        e^{15, 6}_{16, 3},
        e^{15, 7}_{16, 4},
        e^{15, 8}_{16, 5}
        \\
        e^{16, 2}_{14, 6},
        e^{16, 3}_{15, 6},
        e^{16, 4}_{15, 7},
        e^{16, 5}_{15, 8}
      \end{array}
      \right\}
    \end{array}
  \right\},
  \\
  &&
  \hspace{-45pt}
  -
  \Delta
  q^{\frac{5}{2}}
  A_{1}^{\frac{1}{2}}
  A_{2}^{\frac{1}{2}}
  \left\{
    \begin{array}{@{\hspace{-0.8mm}}c@{\hspace{-0.8mm}}}
      e^{ 6, 7}_{12, 2},
      e^{ 6, 8}_{13, 2},
      e^{ 6, 9}_{12, 3},
      e^{ 6,10}_{13, 3},
      e^{ 7, 8}_{14, 2},
      e^{ 7, 9}_{12, 4},
      e^{ 7,11}_{14, 4},
      e^{ 8,10}_{13, 5},
      e^{ 8,11}_{14, 5},
      e^{ 9,10}_{15, 3},
      e^{ 9,11}_{15, 4},
      e^{10,11}_{15, 5}
      \\
      e^{12, 2}_{ 6, 7},
      e^{13, 2}_{ 6, 8},
      e^{12, 3}_{ 6, 9},
      e^{13, 3}_{ 6,10},
      e^{14, 2}_{ 7, 8},
      e^{12, 4}_{ 7, 9},
      e^{14, 4}_{ 7,11},
      e^{13, 5}_{ 8,10},
      e^{14, 5}_{ 8,11},
      e^{15, 3}_{ 9,10},
      e^{15, 4}_{ 9,11},
      e^{15, 5}_{10,11}
    \end{array}
  \right\},
  \\
  &&
  \hspace{-45pt}
  \Delta
  q^{\frac{7}{2}}
  A_{1}^{\frac{1}{2}}
  A_{2}^{\frac{1}{2}}
  \left\{
    \begin{array}{@{\hspace{0mm}}c@{\hspace{0mm}}}
      e^{ 7, 6}_{12, 2},
      e^{ 8, 6}_{13, 2},
      e^{ 8, 7}_{14, 2},
      e^{ 9, 6}_{12, 3},
      e^{ 9, 7}_{12, 4},
      e^{10, 6}_{13, 3},
      e^{10, 8}_{13, 5},
      e^{10, 9}_{15, 3},
      e^{11, 7}_{14, 4},
      e^{11, 8}_{14, 5},
      e^{11, 9}_{15, 4},
      e^{11,10}_{15, 5}
      \\
      e^{12, 2}_{ 7, 6},
      e^{13, 2}_{ 8, 6},
      e^{14, 2}_{ 8, 7},
      e^{12, 3}_{ 9, 6},
      e^{12, 4}_{ 9, 7},
      e^{13, 3}_{10, 6},
      e^{13, 5}_{10, 8},
      e^{15, 3}_{10, 9},
      e^{14, 4}_{11, 7},
      e^{14, 5}_{11, 8},
      e^{15, 4}_{11, 9},
      e^{15, 5}_{11,10}
    \end{array}
  \right\},
  \\
  &&
  \hspace{-45pt}
  \Delta
  q^{4}
  A_{1}^{\frac{1}{2}}
  A_{2}^{\frac{1}{2}}
  \left\{
    -
    \overline{q}^{\frac{1}{2}}
    \left\{
    \begin{array}{@{\hspace{0mm}}c@{\hspace{0mm}}}
      e^{ 6,11}_{14, 3},
      e^{ 7,10}_{13, 4},
      e^{ 8, 9}_{12, 5}
      \\
      e^{14, 3}_{ 6,11},
      e^{13, 4}_{ 7,10},
      e^{12, 5}_{ 8, 9}
    \end{array}
    \right\},
    +
    q^{\frac{1}{2}}
    \left\{
    \begin{array}{@{\hspace{0mm}}c@{\hspace{0mm}}}
      e^{ 9, 8}_{13, 4},
      e^{10, 7}_{12, 5},
      e^{ 6,11}_{15, 2}
      \\
      e^{13, 4}_{ 9, 8},
      e^{12, 5}_{10, 7},
      e^{15, 2}_{ 6,11}
    \end{array}
    \right\}
  \right\},
  \\
  &&
  \hspace{-45pt}
  \Delta
  q^{6}
  A_{1}^{\frac{1}{2}}
  A_{2}^{\frac{1}{2}}
  \left\{
    -
    \overline{q}^{\frac{1}{2}}
    \left\{
    \begin{array}{@{\hspace{0mm}}c@{\hspace{0mm}}}
      e^{ 9, 8}_{14, 3},
      e^{ 7,10}_{15, 2},
      e^{11, 6}_{12, 5}
      \\
      e^{14, 3}_{ 9, 8},
      e^{15, 2}_{ 7,10},
      e^{12, 5}_{11, 6}
    \end{array}
    \right\},
    +
    q^{\frac{1}{2}}
    \left\{
    \begin{array}{@{\hspace{0mm}}c@{\hspace{0mm}}}
      e^{11, 6}_{13, 4},
      e^{10, 7}_{14, 3},
      e^{ 8, 9}_{15, 2}
      \\
      e^{13, 4}_{11, 6},
      e^{14, 3}_{10, 7},
      e^{15, 2}_{ 8, 9}
    \end{array}
    \right\}
  \right\},
  \\
  &&
  \hspace{-45pt}
  \Delta^2
  q^{4}
  A_{1}^{\frac{1}{2}}
  A_{2}^{\frac{1}{2}}
  \left\{
    +
    \overline{q}^{\frac{1}{2}}
    \left\{
    \begin{array}{@{\hspace{0mm}}c@{\hspace{0mm}}}
      e^{ 7,10}_{14, 3},
      e^{ 8, 9}_{13, 4}
      \\
      e^{14, 3}_{ 7,10},
      e^{13, 4}_{ 8, 9}
    \end{array}
    \right\},
    -
    q^{\frac{1}{2}}
    \left\{
    \begin{array}{@{\hspace{0mm}}c@{\hspace{0mm}}}
      e^{ 8, 9}_{14, 3},
      e^{10, 7}_{13, 4}
      \\
      e^{13, 4}_{10, 7},
      e^{14, 3}_{ 8, 9}
    \end{array}
    \right\}
  \right\},
  \\
  &&
  \hspace{-45pt}
  \Delta^2
  q^{\frac{11}{2}}
  (\overline{q}^2 - q^{2})
  A_{1}^{\frac{1}{2}}
  A_{2}^{\frac{1}{2}}
  \left\{
    -
    \overline{q}
    \left\{
    \begin{array}{@{\hspace{0mm}}c@{\hspace{0mm}}}
      e^{ 9, 8}_{15, 2}
      \\
      e^{15, 2}_{ 9, 8}
    \end{array}
    \right\},
    +
    1
    \left\{
    \begin{array}{@{\hspace{0mm}}c@{\hspace{0mm}}}
      e^{10, 7}_{15, 2},
      e^{11, 6}_{14, 3}
      \\
      e^{15, 2}_{10, 7},
      e^{14, 3}_{11, 6}
    \end{array}
    \right\},
    -
    q
    \left\{
    \begin{array}{@{\hspace{0mm}}c@{\hspace{0mm}}}
      e^{11, 6}_{15, 2}
      \\
      e^{15, 2}_{11, 6}
    \end{array}
    \right\}
  \right\},
  \\
  &&
  \hspace{-45pt}
  \Delta^2
  q^{6}
  A_{1}
  A_{2}
  \left\{
    \overline{q}^{2}
    \left\{
    \begin{array}{@{\hspace{0mm}}c@{\hspace{0mm}}}
      e^{12, 5}_{13, 4}
      \\
      e^{13, 4}_{12, 5}
    \end{array}
    \right\},
    -
    \overline{q}
    \left\{
    \begin{array}{@{\hspace{0mm}}c@{\hspace{0mm}}}
      e^{12, 5}_{14, 3}
      \\
      e^{14, 3}_{12, 5}
    \end{array}
    \right\},
    1
    \left\{
    \begin{array}{@{\hspace{0mm}}c@{\hspace{0mm}}}
      e^{12, 5}_{15, 2},
      e^{13, 4}_{14, 3}
      \\
      e^{15, 2}_{12, 5},
      e^{14, 3}_{13, 4}
    \end{array}
    \right\},
    -
    q
    \left\{
    \begin{array}{@{\hspace{0mm}}c@{\hspace{0mm}}}
      e^{13, 4}_{15, 2}
      \\
      e^{15, 2}_{13, 4}
    \end{array}
    \right\},
    q^{2}
    \left\{
    \begin{array}{@{\hspace{0mm}}c@{\hspace{0mm}}}
      e^{14, 3}_{15, 2}
      \\
      e^{15, 2}_{14, 3}
    \end{array}
    \right\}
  \right\},
  \\
  &&
  \hspace{-45pt}
  \Delta^2
  q^{5}
  A_{0}^{\frac{1}{2}}
  A_{1}^{\frac{1}{2}}
  A_{2}^{\frac{1}{2}}
  A_{3}^{\frac{1}{2}}
  \times
  \\
  & &
  \left\{
    \overline{q}^{2}
    \left\{
    \begin{array}{@{\hspace{0mm}}c@{\hspace{0mm}}}
      e^{ 6,11}_{16, 1}
      \\
      e^{16, 1}_{ 6,11}
    \end{array}
    \right\},
    -
    \overline{q}
    \left\{
    \begin{array}{@{\hspace{0mm}}c@{\hspace{0mm}}}
      e^{ 7,10}_{16, 1}
      \\
      e^{16, 1}_{ 7,10}
    \end{array}
    \right\},
    1
    \left\{
    \begin{array}{@{\hspace{0mm}}c@{\hspace{0mm}}}
      e^{ 8, 9}_{16, 1},
      e^{ 9, 8}_{16, 1}
      \\
      e^{16, 1}_{ 8, 9},
      e^{16, 1}_{ 9, 8}
    \end{array}
    \right\},
    -
    q
    \left\{
    \begin{array}{@{\hspace{0mm}}c@{\hspace{0mm}}}
      e^{10, 7}_{16, 1}
      \\
      e^{16, 1}_{10, 7}
    \end{array}
    \right\},
    q^{2}
    \left\{
    \begin{array}{@{\hspace{0mm}}c@{\hspace{0mm}}}
      e^{11, 6}_{16, 1}
      \\
      e^{16, 1}_{11, 6}
    \end{array}
    \right\}
  \right\},
  \\
  &&
  \hspace{-45pt}
  -
  \Delta
  q^{3}
  A_{0}^{\frac{1}{2}}
  A_{3}^{\frac{1}{2}}
  \left\{
    +
    \overline{q}^{\frac{3}{2}}
    \left\{
    \begin{array}{@{\hspace{0mm}}c@{\hspace{0mm}}}
      e^{ 2,15}_{16, 1}
      \\
      e^{16, 1}_{ 2,15}
    \end{array}
    \right\},
    -
    \overline{q}^{\frac{1}{2}}
    \left\{
    \begin{array}{@{\hspace{0mm}}c@{\hspace{0mm}}}
      e^{ 3,14}_{16, 1}
      \\
      e^{16, 1}_{ 3,14}
    \end{array}
    \right\},
    +
    q^{\frac{1}{2}}
    \left\{
    \begin{array}{@{\hspace{0mm}}c@{\hspace{0mm}}}
      e^{ 4,13}_{16, 1}
      \\
      e^{16, 1}_{ 4,13}
    \end{array}
    \right\},
    -
    q^{\frac{3}{2}}
    \left\{
    \begin{array}{@{\hspace{0mm}}c@{\hspace{0mm}}}
      e^{ 5,12}_{16, 1}
      \\
      e^{16, 1}_{ 5,12}
    \end{array}
    \right\}
  \right\},
  \\
  &&
  \hspace{-45pt}
  -
  \Delta^3
  q^{6}
  A_{0}^{\frac{1}{2}}
  A_{1}
  A_{2}
  A_{3}^{\frac{1}{2}}
  \left\{
    +
    \overline{q}^{\frac{3}{2}}
    \left\{
    \begin{array}{@{\hspace{0mm}}c@{\hspace{0mm}}}
      e^{12, 5}_{16, 1}
      \\
      e^{16, 1}_{12, 5}
    \end{array}
    \right\},
    -
    \overline{q}^{\frac{1}{2}}
    \left\{
    \begin{array}{@{\hspace{0mm}}c@{\hspace{0mm}}}
      e^{13, 4}_{16, 1}
      \\
      e^{16, 1}_{13, 4}
    \end{array}
    \right\},
    +
    q^{\frac{1}{2}}
    \left\{
    \begin{array}{@{\hspace{0mm}}c@{\hspace{0mm}}}
      e^{14, 3}_{16, 1}
      \\
      e^{16, 1}_{14, 3}
    \end{array}
    \right\},
    -
    q^{\frac{3}{2}}
    \left\{
    \begin{array}{@{\hspace{0mm}}c@{\hspace{0mm}}}
      e^{15, 2}_{16, 1}
      \\
      e^{16, 1}_{15, 2}
    \end{array}
    \right\}
  \right\},
  \\
  &&
  \hspace{-45pt}
  -
  \Delta
  q^{2 \alpha + \frac{13}{2}}
  A_{2}^{\frac{1}{2}}
  A_{3}^{\frac{1}{2}}
  \left\{
    \begin{array}{@{\hspace{0mm}}c@{\hspace{0mm}}}
      e^{12,13}_{16, 6},
      e^{12,14}_{16, 7},
      e^{13,14}_{16, 8},
      e^{12,15}_{16, 9},
      e^{13,15}_{16,10},
      e^{14,15}_{16,11}
      \\
      e^{16, 6}_{12,13},
      e^{16, 7}_{12,14},
      e^{16, 8}_{13,14},
      e^{16, 9}_{12,15},
      e^{16,10}_{13,15},
      e^{16,11}_{14,15}
    \end{array}
  \right\},
  \\
  &&
  \hspace{-45pt}
  \Delta
  q^{2 \alpha + \frac{15}{2}}
  A_{2}^{\frac{1}{2}}
  A_{3}^{\frac{1}{2}}
  \left\{
    \begin{array}{@{\hspace{0mm}}c@{\hspace{0mm}}}
      e^{13,12}_{16, 6},
      e^{14,12}_{16, 7},
      e^{14,13}_{16, 8},
      e^{15,12}_{16, 9},
      e^{15,13}_{16,10},
      e^{15,14}_{16,11}
      \\
      e^{16, 6}_{13,12},
      e^{16, 7}_{14,12},
      e^{16, 8}_{14,13},
      e^{16, 9}_{15,12},
      e^{16,10}_{15,13},
      e^{16,11}_{15,14}
    \end{array}
  \right\},
\end{eqnarray*}
and the left handle $C$ has $16$ diagonal components:
\begin{eqnarray*}
  C
  =
  q^{- 4 \alpha - 6}
  \left\{
    \begin{array}{l}
      +
      q^{6}
      \left\{
                         e^{1}_{1}
      \right\}
      \\
      -
      q^{3}
      \left\{
        q^{3}            e^{2}_{2}, \; \;
        q                e^{3}_{3}, \; \;
        \overline{q}     e^{4}_{4}, \; \;
        \overline{q}^{3} e^{5}_{5}
      \right\}
      \\
      +
      \left\{
        q^{4}            e^{6}_{6}, \; \;
        q^{2}            e^{7}_{7}, \; \;
                         e^{8}_{8}, \; \;
                         e^{9}_{9}, \; \;
        \overline{q}^{2} e^{10}_{10}, \; \;
        \overline{q}^{4} e^{11}_{11}
      \right\}
      \\
      -
      \overline{q}^{3}
      \left\{
        \overline{q}^{3} e^{15}_{15}, \; \;
        \overline{q}     e^{14}_{14}, \; \;
        q                e^{13}_{13}, \; \;
        q^{3}            e^{12}_{12}
      \right\}
      \\
      +
      \overline{q}^{6}
      \left\{
                       e^{16}_{16}
      \right\}
   \end{array}
  \right\},
\end{eqnarray*}
using the scaling factors:
\begin{eqnarray*}
  \kappa_{\sigma}
  =
  q^{- 4 \alpha},
  \qquad \qquad
  \kappa_{C}
  =
  q^{- 4 \alpha}.
\end{eqnarray*}

\vspace{5mm}


\section{Evaluations of $LG^3$ and $LG^4$}
\label{app:LG34}

Below, we present evaluations for $LG^3$ and $LG^{4}$ for a few $2$ and
$3$-braid links.  We use the same naming conventions for links as those
of \cite{DeWit:99a}, except that we denote by $2^2_{1a}$ and $2^2_{1b}$
the $2$ component links determined respectively by the braids
$\sigma_1^{\pm 1}$.

To present evaluations of $LG^m$,
we use a similar convention to that of \cite{DeWit:99a}.
The expression $s_0 (A_0(q)), s_1 (A_1(q)), \dots, s_r (A_r(q))$,
where the $s_i$ are signs $\pm$ and the $A_i(q)$ are
integer-coefficient Laurent polynomials in $q$, is intended to be read:
\begin{eqnarray*}
    s_0 (A_0(q))
  + s_1 (\overline{P}^1 + P^1) (A_1(q))
  + \dots
  + s_r (\overline{P}^r + P^r) (A_r(q)).
\end{eqnarray*}
In these expressions $(A_i(q))$ is only a list of terms of
$A_i(q)$ rather than an explicit sum, viz we have written
$(x_1, x_2, \dots, x_s)$ for $(x_1 + x_2 + \cdots + x_s)$. Recall that,
for $LG^{m,n}$ (fixing $m,n$), we are using the variable
$P=p^2=q^{2\alpha+m-n}$; viz for
$LG^3\equiv LG^{3,1}$ we use $P=q^{2\alpha+2}$, and for
$LG^4\equiv LG^{4,1}$ we use $P=q^{2\alpha+3}$.
For multicomponent links, if the polynomial is not invariant under
$P \mapsto \overline{P}$, we write it out in full.
This situation only occurs here for $LG^{3}$ for links of $2$
components.

Within the $q$-polynomials, the same general behaviours of the
coefficients as reported for $LG^{2}$ in \cite{DeWit:99a} are seen.
These calculations were performed on \texttt{SUN Ultra 60} UNIX based
workstations, with a main memory of 256Mb, and the larger calculations
sometimes used all of this memory.

\subsection*{Evaluations of $LG^{3}$}

\small

\begin{eqnarray*}
  \hspace{-30pt}
  \begin{array}{rcl}
    LG^3_{3_{1}}
    & {\hspace{-10pt} = \hspace{-10pt}} &
    -
    (q^{2}, 2 q^{4}, 3 q^{6}, q^{8}),
    +
    (q^{2}, 2 q^{4}, 2 q^{6}, q^{8}),
    -
    (q^{4}, q^{6}, q^{8}),
    +
    (q^{6})
    \\
    \\
    LG^3_{4_1}
    & {\hspace{-10pt} = \hspace{-10pt}} &
    +
    (5 \overline{q}^{4}, 9 \overline{q}^{2}, 17, 9 q^{2}, 5 q^{4}),
    -
    (2 \overline{q}^{4}, 8 \overline{q}^{2}, 10, 8 q^{2}, 2 q^{4}),
    +
    (3 \overline{q}^{2}, 3, 3 q^{2}),
    -
    (1)
    \\
    \\
    LG^3_{5_1}
    & {\hspace{-10pt} = \hspace{-10pt}} &
    +
    (q^{4}, 4 q^{6}, 5 q^{8}, 5 q^{10}, 3 q^{12}, q^{14}),
    -
    (q^{4}, 3 q^{6}, 5 q^{8}, 5 q^{10}, 3 q^{12}, q^{14}),
    \\
    & &
    +
    (q^{4}, 2 q^{6}, 4 q^{8}, 4 q^{10}, 3 q^{12}, q^{14}),
    -
    (q^{6}, 2 q^{8}, 3 q^{10}, 3 q^{12}, q^{14}),
    \\
    & &
    +
    (q^{8}, 2 q^{10}, 2 q^{12}, q^{14}),
    -
    (q^{10}, q^{12}, q^{14}),
    +
    (q^{12})
    \\
    \\
    LG^3_{5_2}
    & {\hspace{-10pt} = \hspace{-10pt}} &
    -
    (7 q^{2}, 17 q^{4}, 32 q^{6}, 25 q^{8}, 15 q^{10}, 3 q^{12}),
    +
    (4 q^{2}, 15 q^{4}, 23 q^{6}, 22 q^{8}, 11 q^{10}, 3 q^{12}),
    \\
    & &
    -
    (7 q^{4}, 10 q^{6}, 12 q^{8}, 5 q^{10}, 2 q^{12}),
    +
    (4 q^{6}, 2 q^{8}, 2 q^{10})
  \end{array}
\end{eqnarray*}

\begin{eqnarray*}
  \hspace{-30pt}
  \begin{array}{rcl}
    LG^3_{6_2}
    & {\hspace{-10pt} = \hspace{-10pt}} &
    -
    (7 \overline{q}^{2}, 29, 60 q^{2}, 74 q^{4}, 60 q^{6}, 26 q^{8},
      5 q^{10}),
    \\
    & &
    +
    (6 \overline{q}^{2}, 24, 52 q^{2}, 67 q^{4}, 54 q^{6}, 26 q^{8},
      5 q^{10}),
    \\
    & &
    -
    (2 \overline{q}^{2}, 14, 32 q^{2}, 48 q^{4}, 41 q^{6}, 23 q^{8},
      5 q^{10}),
    +
    (4, 15 q^{2}, 24 q^{4}, 28 q^{6}, 14 q^{8}, 5 q^{10}),
    \\
    & &
    -
    (4 q^{2}, 10 q^{4}, 12 q^{6}, 8 q^{8}, 2 q^{10}),
    +
    (3 q^{4}, 3 q^{6}, 3 q^{8}),
    -
    (q^{6})
    \\
    \\
    LG^3_{6_3}
    & {\hspace{-10pt} = \hspace{-10pt}} &
    +
    (9 \overline{q}^{6}, 52 \overline{q}^{4}, 106 \overline{q}^{2}, 145,
     106 q^{2}, 52 q^{4}, 9 q^{6}),
    \\
    & &
    -
    (9 \overline{q}^{6}, 42 \overline{q}^{4}, 96 \overline{q}^{2}, 120,
     96 q^{2}, 42 q^{4}, 9 q^{6}),
    \\
    & &
    +
    (6 \overline{q}^{6}, 26 \overline{q}^{4}, 63 \overline{q}^{2}, 77,
     63 q^{2}, 26 q^{4}, 6 q^{6}),
    \\
    & &
    -
    (2 \overline{q}^{6}, 13 \overline{q}^{4}, 28 \overline{q}^{2},
     40, 28 q^{2}, 13 q^{4}, 2 q^{6}),
    +
    (4 \overline{q}^{4}, 10 \overline{q}^{2}, 14, 10 q^{2}, 4 q^{4}),
    \\
    & &
    -
    (3 \overline{q}^{2}, 3, 3 q^{2}),
    +
    (1)
  \end{array}
\end{eqnarray*}

\begin{eqnarray*}
  \hspace{-30pt}
  \begin{array}{rcl}
    LG^3_{2^2_{1a}}
    & {\hspace{-10pt} = \hspace{-10pt}} &
    - (\overline{p}^{1} - p^{1})
    (q, q^3, q^5),
    + (\overline{p}^{3} - p^{3})
    (q^3)
    \\
    \\
    LG^3_{2^2_{1b}}
    & {\hspace{-10pt} = \hspace{-10pt}} &
    + (\overline{p}^{1} - p^{1})
    (\overline{q}^5, \overline{q}^3, \overline{q}^1),
    - (\overline{p}^{3} - p^{3})
    (\overline{q}^3)
    \\
    \\
    LG^3_{4^2_{1a}}
    & {\hspace{-10pt} = \hspace{-10pt}} &
    + (\overline{p}^{1} - p^{1})
    (q^{3}, 3 q^{5}, 4 q^{7}, 3 q^{9}, q^{11}),
    - (\overline{p}^{3} - p^{3})
    (q^{3}, 2 q^{5}, 3 q^{7}, 3 q^{9}, q^{11})
    \\
    & &
    + (\overline{p}^{5} - p^{5})
    (q^{5}, 2 q^{7}, 2 q^{9}, q^{11}),
    - (\overline{p}^{7} - p^{7})
    (q^{7}, q^{9}, q^{11}),
    + (\overline{p}^{9} - p^{9})
    (q^{9})
    \\
    \\
    LG^3_{4^2_{1b}}
    & {\hspace{-10pt} = \hspace{-10pt}} &
    - (\overline{p}^{1} - p^{1})
    (4 q^{1}, 6 q^{3}, 8 q^{5}, 4 q^{7}, 2 q^{9}),
    + (\overline{p}^{3} - p^{3})
    (4 q^{3}, 2 q^{5}, 2 q^{7})
    \\
    \\
    LG^3_{5^2_{1}}
    & {\hspace{-10pt} = \hspace{-10pt}} &
    + (\overline{p}^{1} - p^{1})
    (5 \overline{q}^{7}, 23 \overline{q}^{5}, 36 \overline{q}^{3},
      37 \overline{q}^{1}, 19 q^{1}, 6 q^{3}),
    \\
    & &
    - (\overline{p}^{3} - p^{3})
    (5 \overline{q}^{7}, 14 \overline{q}^{5}, 28 \overline{q}^{3},
      22 \overline{q}^{1}, 13 q^{1}, 2 q^{3})
    \\
    & &
    + (\overline{p}^{5} - p^{5})
    (2 \overline{q}^{7}, 8 \overline{q}^{5}, 12 \overline{q}^{3},
      10 \overline{q}^{1}, 4 q^{1})
    \\
    & &
    - (\overline{p}^{7} - p^{7})
    (3 \overline{q}^{5}, 3 \overline{q}^{3}, 3 \overline{q}^{1}),
    + (\overline{p}^{9} - p^{9})
    (\overline{q}^{3})
    \\
    \\
    LG^3_{6^2_{1}}
    & {\hspace{-10pt} = \hspace{-10pt}} &
    - (\overline{p}^{1} - p^{1})
    (q^{5}, 4 q^{7}, 6 q^{9}, 7 q^{11}, 5 q^{13}, 3 q^{15}, q^{17})
    \\
    & &
    + (\overline{p}^{3} - p^{3})
    (q^{5}, 3 q^{7}, 6 q^{9}, 6 q^{11}, 5 q^{13}, 3 q^{15}, q^{17})
    \\
    & &
    - (\overline{p}^{5} - p^{5})
    (q^{5}, 2 q^{7}, 4 q^{9}, 5 q^{11}, 5 q^{13}, 3 q^{15}, q^{17})
    \\
    & &
    + (\overline{p}^{7} - p^{7})
    (q^{7}, 2 q^{9}, 4 q^{11}, 4 q^{13}, 3 q^{15}, q^{17})
    \\
    & &
    - (\overline{p}^{9} - p^{9})
    (q^{9}, 2 q^{11}, 3 q^{13}, 3 q^{15}, q^{17}),
    + (\overline{p}^{11} - p^{11})
    (q^{11}, 2 q^{13}, 2 q^{15}, q^{17})
    \\
    & &
    - (\overline{p}^{13} - p^{13})
    (q^{13}, q^{15}, q^{17}),
    + (\overline{p}^{15} - p^{15})
    (q^{15} )
    \\
    \\
    LG^3_{6^2_{2}}
    & {\hspace{-10pt} = \hspace{-10pt}} &
    + (\overline{p}^{1} - p^{1})
    (7 q^{3}, 26 q^{5}, 52 q^{7}, 59 q^{9}, 43 q^{11}, 17 q^{13}, 3 q^{15})
    \\
    & &
    - (\overline{p}^{3} - p^{3})
    (4 q^{3}, 18 q^{5}, 36 q^{7}, 48 q^{9}, 34 q^{11}, 16 q^{13}, 3 q^{15})
    \\
    & &
    + (\overline{p}^{5} - p^{5})
    (7 q^{5}, 19 q^{7}, 27 q^{9}, 23 q^{11}, 11 q^{13}, 3 q^{15})
    \\
    & &
    - (\overline{p}^{7} - p^{7})
    (7 q^{7}, 10 q^{9}, 12 q^{11}, 5 q^{13}, 2 q^{15}),
    + (\overline{p}^{9} - p^{9})
    (4 q^{9}, 2 q^{11}, 2 q^{13})
    \\
    \\
    \\
    LG^3_{6^3_{1}}
    & {\hspace{-10pt} = \hspace{-10pt}} &
    -
    (10 \overline{q}^{2}, 38, 82 q^{2}, 100 q^{4}, 84 q^{6},
     36 q^{8}, 8 q^{10}),
    \\
    & &
    +
    (7 \overline{q}^{2}, 32, 69 q^{2}, 91 q^{4}, 73 q^{6}, 36 q^{8},
     7 q^{10}),
    \\
    & &
    -
    (2 \overline{q}^{2}, 18, 42 q^{2}, 63 q^{4}, 53 q^{6}, 29 q^{8},
     6 q^{10}),
    +
    (5, 19 q^{2}, 30 q^{4}, 33 q^{6}, 16 q^{8}, 5 q^{10}),
    \\
    & &
    -
    (5 q^{2}, 11 q^{4}, 13 q^{6}, 8 q^{8}, 2 q^{10}),
    +
    (3 q^{4}, 3 q^{6}, 3 q^{8}),
    -
    (q^{6})
    \\
    \\
    LG^3_{6^3_{2}}
    & {\hspace{-10pt} = \hspace{-10pt}} &
    +
    (18 \overline{q}^{6}, 100 \overline{q}^{4}, 206 \overline{q}^{2},
     276, 206 q^{2}, 100 q^{4}, 18 q^{6}),
    \\
    & &
    -
    (17 \overline{q}^{6}, 81 \overline{q}^{4}, 183 \overline{q}^{2},
     230, 183 q^{2}, 81 q^{4}, 17 q^{6}),
    \\
    & &
    +
    (11 \overline{q}^{6}, 48 \overline{q}^{4}, 117 \overline{q}^{2},
     143, 117 q^{2}, 48 q^{4}, 11 q^{6}),
    \\
    & &
    -
    (3 \overline{q}^{6}, 23 \overline{q}^{4}, 49 \overline{q}^{2},
     70, 49 q^{2}, 23 q^{4}, 3 q^{6}),
    +
    (6 \overline{q}^{4}, 16 \overline{q}^{2}, 22, 16 q^{2}, 6 q^{4}),
    \\
    & &
    -
    (4 \overline{q}^{2}, 4, 4 q^{2}),
    +
    (1)
      \end{array}
\end{eqnarray*}

\begin{eqnarray*}
  \hspace{-30pt}
  \begin{array}{rcl}
    LG^3_{6^3_{3}}
    & {\hspace{-10pt} = \hspace{-10pt}} &
    +
    (2 q^{6}, 2 q^{8}, 2 q^{10}),
    -
    (q^{4}, 2 q^{6}, 3 q^{8}, 2 q^{10}, q^{12}),
    +
    (q^{4}, 2 q^{6}, 3 q^{8}, 2 q^{10}, q^{12}),
    \\
    & &
    -
    (q^{6}, q^{8}, q^{10}, q^{12}),
    +
    (q^{10}, q^{12}, q^{14}),
    -
    (q^{10}, q^{12}, q^{14}),
    +
    (q^{12})
  \end{array}
\end{eqnarray*}

\normalsize

\subsection*{Evaluations of $LG^{4}$}

\small

\begin{eqnarray*}
  \hspace{-30pt}
  \begin{array}{rcl}
    \hspace{-30pt}
    LG^4_{2^2_{1a}}
    & {\hspace{-10pt} = \hspace{-10pt}} &
    +
    (q^{2}, q^{4}, 2 q^{6}, q^{8}, q^{10}),
    -
    (q^{3}, q^{5}, q^{7}, q^{9}),
    +
    (q^{6})
    \\
    \\
    LG^4_{2^2_{1b}}
    & {\hspace{-10pt} = \hspace{-10pt}} &
    +
    (\overline{q}^{10}, \overline{q}^{8}, 2 \overline{q}^{6}, \overline{q}^{4},
     \overline{q}^{2}),
    -
    (\overline{q}^{9}, \overline{q}^{7}, \overline{q}^{5}, \overline{q}^{3}),
    +
    (\overline{q}^{6})
    \\
    \\
    LG^4_{3_{1}}
    & {\hspace{-10pt} = \hspace{-10pt}} &
    +
    (q^{4}, 2 q^{6}, 4 q^{8}, 4 q^{10}, 5 q^{12}, 2 q^{14},
     q^{16}),
    -
    (q^{5}, 2 q^{7}, 4 q^{9}, 4 q^{11}, 3 q^{13}, 2 q^{15}),
    \\
    & &
    +
    (  q^{6}, 2 q^{8}, 2 q^{10}, 3 q^{12}, q^{14}, q^{16}),
    -
    (q^{9}, q^{11}, q^{13}, q^{15}),
    +
    (  q^{12})
    \\
    \\
    LG^4_{4^2_{1a}}
    & {\hspace{-10pt} = \hspace{-10pt}} &
    +
    (q^{6}, 2 q^{8}, 5 q^{10}, 8 q^{12}, 10 q^{14}, 8 q^{16},
     7 q^{18}, 2 q^{20}, q^{22}),
    \\
    & &
    -
    (q^{7}, 3 q^{9}, 6 q^{11}, 8 q^{13}, 9 q^{15}, 7 q^{17},
     4 q^{19}, 2 q^{21}),
    \\
    & &
    +
    (q^{8}, 3 q^{10}, 5 q^{12}, 7 q^{14}, 6 q^{16}, 6 q^{18},
     2 q^{20}, q^{22}), 
    \\
    & &
    -
    (q^{9}, 2 q^{11}, 3 q^{13}, 5 q^{15}, 4 q^{17}, 3 q^{19},
     2 q^{21}),
    +
    (q^{12}, 2 q^{14}, 2 q^{16}, 3 q^{18}, q^{20}, q^{22}),
    \\
    & &
    -
    (q^{15}, q^{17}, q^{19}, q^{21}), 
    +
    (q^{18})
    \\
    \\
    LG^4_{5_{1}}
    & {\hspace{-10pt} = \hspace{-10pt}} &
    +
    (q^8, 2 q^{10}, 6 q^{12}, 10 q^{14}, 15 q^{16}, 16 q^{18},
     15 q^{20}, 10 q^{22}, 7 q^{24}, 2 q^{26}, q^{28}),
    \\
    & &
    -
    (q^9, 3 q^{11}, 7 q^{13}, 12 q^{15}, 15 q^{17}, 15 q^{19},
     13 q^{21}, 8 q^{23}, 4 q^{25}, 2 q^{27}),
    \\
    & &
    +
    (q^{10}, 4 q^{12}, 7 q^{14}, 11 q^{16}, 13 q^{18}, 13 q^{20},
     9 q^{22}, 7 q^{24}, 2 q^{26}, q^{28}),
    \\
    & &
    -
    (q^{11}, 3 q^{13}, 6 q^{15}, 9 q^{17}, 10 q^{19}, 10 q^{21},
     7 q^{23}, 4 q^{25}, 2 q^{27}),
    \\
    & &
    +
    (q^{12}, 2 q^{14}, 4 q^{16}, 6 q^{18}, 7 q^{20}, 6 q^{22}, 6 q^{24},
     2 q^{26}, q^{28}),
    \\
    & &
    -
    (q^{15}, 2 q^{17}, 3 q^{19}, 5 q^{21}, 4 q^{23}, 3 q^{25},
     2 q^{27}),
    \\
    & &
    +
    (q^{18}, 2 q^{20}, 2 q^{22}, 3 q^{24}, q^{26}, q^{28}),
    -
    (q^{21}, q^{23}, q^{25}, q^{27}),
    +
    (q^{24})
  \end{array}
\end{eqnarray*}

\normalsize


\bibliographystyle{plain}
\bibliography{DeWit99e}

\begin{thebibliography}{10}

\bibitem{AkutsuDeguchiOhtsuki:92}
Yasuhiro Akutsu, Tetsuo Deguchi, and Tomotada Ohtsuki.
\newblock Invariants of colored links.
\newblock {\em Journal of Knot Theory and its Ramifications}, 1(2):161--184,
  1992.

\bibitem{DeWit:98}
David {De Wit}.
\newblock Explicit construction of the representation of the braid generator
  {$\sigma$} associated with the one-parameter family of minimal typical
  highest weight {$(0,0\,|\,\alpha)$} representations of {$U_q[gl(2|1)]$} and
  its use in the evaluation of the {L}inks--{G}ould two-variable {L}aurent
  polynomial invariant of oriented $(1,1)$ tangles, 25 November 1998.
\newblock PhD thesis, Department of Mathematics, The University of Queensland,
  Australia. \texttt{math/9909063}.

\bibitem{DeWit:99c}
David {De Wit}.
\newblock Automatic construction of explicit {R} matrices for the one-parameter
  family of irreducible typical highest weight {$(\dot{0}_m|\alpha)$}
  representations of {$U_q[gl(m|1)]$}.
\newblock In preparation. Methods leading to results in \cite{DeWit:99d}, 2000.

\bibitem{DeWit:99a}
David {De Wit}.
\newblock Automatic evaluation of the {L}inks--{G}ould invariant for all prime
  knots of up to $10$ crossings.
\newblock {\em Journal of Knot Theory and its Ramifications}, 9(3):311--339,
  May 2000.
\newblock RIMS-1235, \texttt{math/9906059}.

\bibitem{DeWit:99d}
David {De Wit}.
\newblock Four easy pieces -- explicit {R} matrices from the
  {$(\dot{0}_m|\alpha)$} highest weight representations of {$U_q[gl(m|1)]$}.
\newblock Results of the methods in \cite{DeWit:99c}. Under consideration.
  \texttt{math/0005049}, 5 May 2000.

\bibitem{DeWitKauffmanLinks:99a}
David {De Wit}, Louis~H Kauffman, and Jon~R Links.
\newblock On the {L}inks--{G}ould invariant of links.
\newblock {\em Journal of Knot Theory and its Ramifications}, 8(2):165--199,
  March 1999.
\newblock \texttt{math/9811128}.

\bibitem{DeliusGouldLinksZhang:95b}
Gustav~W Delius, Mark~D Gould, Jon~R Links, and Yao-Zhong Zhang.
\newblock Solutions of the {Y}ang--{B}axter equation with extra non-additive
  parameters {II}: {$ U_q ( gl ( m | n ) ) $}.
\newblock {\em Journal of Physics A. Mathematical and General},
  28(21):6203--6210, 1995.

\bibitem{Hennings:91}
M~A Hennings.
\newblock {H}opf algebras and regular isotopy invariants for link diagrams.
\newblock {\em Mathematical Proceedings of the Cambridge Philosophical
  Society}, 109(1):59--77, January 1991.

\bibitem{HosteThistlethwaiteWeeks:98}
Jim Hoste, Morwen Thistlethwaite, and Jeff Weeks.
\newblock The first $1,701,936$ knots.
\newblock {\em The Mathematical Intelligencer}, 20(4):33--48, 1998.

\bibitem{Kauffman:87b}
Louis~H Kauffman.
\newblock State models and the {J}ones polynomial.
\newblock {\em Topology}, 26(3):395--407, 1987.

\bibitem{Kauffman:93}
Louis~H Kauffman.
\newblock {\em {Knots and Physics}}.
\newblock World Scientific, Singapore, 2nd edition, 1993.

\bibitem{Kauffman:97a}
Louis~H Kauffman.
\newblock Knots and diagrams.
\newblock In Shin'ichi Suzuki, editor, {\em Lectures at {K}nots96}, pages
  123--194. World Scientific, 1997.

\bibitem{KhoroshkinTolstoy:91}
Sergei Khoroshkin and Valerij~N Tolstoy.
\newblock Universal ${R}$-matrix for quantized (super)algebras.
\newblock {\em Communications in Mathematical Physics}, 141(3):599--617, 1991.

\bibitem{KirbyMelvin:91}
Robion Kirby and Paul Melvin.
\newblock The $3$-manifold invariants of {W}itten and {R}eshetikhin-{T}uraev
  for {$\mathrm{sl}(2, {\mathbb{C}})$}.
\newblock {\em Inventiones Mathematicae}, 105(3):473--545, 1991.

\bibitem{LiaoSong:91}
Li~Liao and Xing~Chang Song.
\newblock Quantum {L}ie superalgebras and ``nonstandard'' braid group
  representations.
\newblock {\em Modern Physics Letters A}, 6(11):959--968, 1991.

\bibitem{Lickorish:97}
W~B~Raymond Lickorish.
\newblock {\em An Introduction to Knot Theory}.
\newblock Springer-Verlag, New York, 1997.

\bibitem{LinksGould:92b}
Jon~R Links and Mark~D Gould.
\newblock Two variable link polynomials from quantum supergroups.
\newblock {\em Letters in Mathematical Physics}, 26(3):187--198, November 1992.

\bibitem{LinksGould:96b}
Jon~R Links and Mark~D Gould.
\newblock A {$ q $}-superdimension formula for irreps of type {I} quantum
  superalgebras.
\newblock {\em Journal of Mathematical Physics}, 37(1):484--492, January 1996.

\bibitem{LinksGouldZhang:93}
Jon~R Links, Mark~D Gould, and Rui~Bin Zhang.
\newblock Quantum supergroups, link polynomials and representation of the braid
  generator.
\newblock {\em Reviews in Mathematical Physics}, 5(2):345--361, 1993.

\bibitem{Matveev:87}
S~V Matveev.
\newblock Generalized surgeries of three-dimensional manifolds and
  representations of homology spheres.
\newblock {\em Akademiya Nauk Soyuza SSR. Matematicheskie Zametki},
  42(2):268--278, 345, 1987.
\newblock In Russian. English translation: Mathematical Notes 42 (1987), no.
  1-2, 651--656.

\bibitem{MortonCromwell:96}
Hugh~R Morton and Peter~R Cromwell.
\newblock Distinguishing mutants by knot polynomials.
\newblock {\em Journal of Knot Theory and its Ramifications}, 5(2):225--238,
  1996.

\bibitem{MurakamiNakanishi:89}
Hitoshi Murakami and Yasutaka Nakanishi.
\newblock On a certain move generating link-homology.
\newblock {\em Mathematische Annalen}, 284(1):75--89, 1989.

\bibitem{ReshetikhinTuraev:90}
Nikolai~Yu Reshetikhin and Vladimir~G Turaev.
\newblock Ribbon graphs and their invariants derived from quantum groups.
\newblock {\em Communications in Mathematical Physics}, 127:1--26, 1990.

\bibitem{RozanskySaleur:92}
Lev Rozansky and Hubert Saleur.
\newblock Quantum field theory for the multi-variable {A}lexander-{C}onway
  polynomial.
\newblock {\em Nuclear Physics. B}, 376(3):461--509, 1992.

\bibitem{Zhang:93}
Rui~Bin Zhang.
\newblock Finite dimensional irreducible representations of the quantum
  supergroup {$U_q(gl(m|n))$}.
\newblock {\em Journal of Mathematical Physics}, 34(3):1236--1254, March 1993.

\end{thebibliography}

\end{document}